\documentclass[reqno]{amsart}
\usepackage{amsmath,amsfonts,amssymb,epsf,graphicx}
\usepackage{float,enumerate}

\newcommand{\Fp}{{\lower.10em\hbox{${\mathbb F}$}}\kern-.1em\hbox{}_p}

\newcommand \Sonic{\Gamma_{sonic}}

\newcommand{\PtUpL}{{P_1}}
\newcommand{\PtLwL}{{P_2}}
\newcommand{\PtLwR}{{P_3}}
\newcommand{\PtUpR}{{P_4}}

\newcommand{\Csect}{{\Lambda}}
\newcommand{\Noz}{{\Omega}}
\newcommand{\CylToNoz}{{\Psi}}

\newcommand{\Vtan}{{\bf V}'}
\newcommand{\Vpl}{{\bf V}^-}
\newcommand{\Vmn}{{\bf V}^+}
\newcommand{\Vpm}{V^\pm}
\newcommand{\Gm}{\Gamma}
\newcommand{\e}{{\varepsilon}}

\renewcommand{\r}{\rho}

\newcommand{\D}{{\mathcal D}}
\newcommand{\R}{{\mathbb R}} 
\newcommand{\M}{{\mathcal M}}
\newcommand{\bR}{{\R}} 
\renewcommand{\u}{{\bf u}}
\renewcommand{\v}{{\bf v}}
\renewcommand{\a}{\alpha}

\newcommand{\del}{\partial}
\newcommand{\be}{\begin{equation}}
\newcommand{\ea}{\end{equation}}
\newcommand{\bea}{\begin{eqnarray}}
\newcommand{\eea}{\end{eqnarray}}

\newcommand{\DM}{{\mathcal DM}}

\newcommand{\uu}{{\bf u}}
\newcommand{\mm}{{\bf m}}
\newcommand{\vv}{{\bf v}}

\newcommand{\aaa}{{\bf a}}
\newcommand{\bb}{{\bf b}}
\newcommand{\ww}{{\bf w}}
\newcommand{\xx}{{\bf x}}

\newcommand{\ff}{{\bf f}}

\newcommand{\pp}{{\bf p}}

\newcommand{\rr}{{\bf r}}
\newcommand{\qq}{{\bf q}}

\newcommand{\FF}{{\bf F}}
\newcommand{\AAA}{{\bf A}}
\newcommand{\BBB}{{\bf B}}
\newcommand{\CCC}{{\bf C}}
\newcommand{\OO}{{\bf O}}
\newcommand{\GG}{{\bf G}}
\newcommand{\bS}{{\mathcal S}}
\newcommand{\BS}{{\bf S}}

\newcommand{\QQ}{{\bf Q}}

\newcommand{\VV}{{\bf V}}

\newcommand{\divg}{ \mbox{div\,}}

\newcommand{\xpr}{{\bf x}'}  
\newcommand{\defd}{:=}

\newcommand{\epsP}{\sigma}

\renewcommand{\d}{\partial}
\newcommand{\n}{{\bf n}}

\newcommand{\bnu}{{\bf n}}
\newtheorem{theorem}{Theorem}[section]

\numberwithin{equation}{section}

\begin{document}
\title[Multidimensional Conservation Laws]{Multidimensional Conservation Laws: \\Overview, Problems, and Perspective}
\author{GUI-QIANG G. Chen}
\address{Gui-Qiang G. Chen, Oxford Centre for Nonlinear PDE, Mathematical Institute, University of Oxford,
         Oxford, OX1 3LB, UK; and Department of Mathematics, Northwestern University,
         Evanston, IL 60208, USA}
\email{\tt chengq@maths.ox.ac.uk}

\begin{abstract}
Some of recent important developments are overviewed, several longstanding
open problems are discussed, and a perspective is presented for the
mathematical theory of multidimensional conservation laws.
Some basic features and phenomena of multidimensional
hyperbolic conservation laws are revealed,
and some samples of multidimensional systems/models and related important
problems are presented and analyzed with emphasis on the prototypes that have been
solved or may be expected to be solved rigorously at least for some cases.
In particular, multidimensional steady supersonic problems and transonic problems,
shock reflection-diffraction problems, and related effective nonlinear approaches
are analyzed. A theory of divergence-measure vector fields and related analytical
frameworks for the analysis of entropy solutions are discussed.
\end{abstract}


\maketitle

\bigskip
\bigskip
\bigskip

\tableofcontents

\bigskip
\section{Introduction}\label{s1}
We overview some of recent important developments, discuss several longstanding
open problems, and present a perspective for the mathematical theory of
multidimensional (M-D, for short) conservation laws.

{\it Hyperbolic Conservation Laws}, quasilinear hyperbolic systems
in divergence form, are one of the most important classes of
nonlinear partial differential equations (PDEs for short).
They typically take the form:
\begin{equation}
\del_t \uu +\nabla_\xx \cdot \ff(\uu)=0,
\qquad \uu\in\R^m,
\label{cons}
\end{equation}
for $(t, \xx)\in \R_+^{d+1}:=\R_+\times \R^d:=[0, \infty)\times\R^d$,
where
$\nabla_\xx=(\partial_{x_1},\dots, \partial_{x_d})$, and
$
\ff=(\ff_1,\dots,\ff_d):\,\, \R^m\to (\R^m)^d
$
is a nonlinear mapping with
$\ff_i: \R^m\to\R^m$ for $i=1,\dots, d$. Another prototypical form is
\begin{equation}\label{cons-2}
\partial_t A_0(u_t, \nabla_\xx u)+\nabla_\xx\cdot \AAA(u_t, \nabla_\xx u)=0,
\end{equation}
where $A_0(p_0, \pp):\R\times \R^{d}\to \R$ and $\AAA(p_0, \pp):\R\times \R^{d}\to \R^d$
with $\pp=(p_1,..., p_d)\in\R^d$ are $C^1$ nonlinear mappings.

The hyperbolicity for (\ref{cons}) requires that, for any $\n\in \bS^{d-1}$,
\begin{eqnarray}\label{hyper:1}
&&(\nabla_\ww \ff(\ww) \cdot {\n})_{m\times m}
\mbox{have $m$ real eigenvalues $\lambda_i(\ww; \n)$,
$1\le i\le m$.}
\end{eqnarray}
We say that system (\ref{cons})
is hyperbolic in a state domain $\D$ if condition (\ref{hyper:1})
holds for any $\ww\in \D$ and $\n\in \bS^{d-1}$.

The hyperbolicity for (\ref{cons-2}) requires that, rewriting (\ref{cons-2})
into second-order nondivergence form such that the terms
$u_{x_ix_j}$ have corresponding coefficients $a_{ij}, i,j=0,1,...,d$,
all the eigenvalues of the matrix $(a_{ij})_{d\times d}$ are real, and
one eigenvalue has a different sign from the other $d$ eigenvalues.

\smallskip
An archetype of nonlinear hyperbolic systems of conservation laws
is the Euler equations for compressible fluids
in $\R^d$, which are a system of $d+2$ conservation laws:
\begin{equation}
\left\{\begin{array}{l}
\del_t\rho +\nabla_\xx\cdot \mm=0,\\
\del_t\mm+\nabla_\xx\cdot\big(\frac{\mm\otimes\mm}{\rho}\big)
  +\nabla_\xx p=0,\\
\del_t (\rho E)+\nabla_\xx\cdot\big(\mm(E+ p/\rho)\big)=0.
\end{array} \right.
\label{Euler1}
\end{equation}
System (\ref{Euler1}) is closed by the constitutive relations:
\begin{equation}
p=p(\rho, e), \qquad
E=\frac{1}{2}\frac{|\mm|^2}{\rho^2}+e.
\label{equilibrium}
\end{equation}
In (\ref{Euler1})--(\ref{equilibrium}), $\tau={1}/{\rho}$ is the
deformation gradient (the specific volume for fluids, the strain for solids),
$\mm=(m_1,\dots,m_d)^\top$ is the fluid momentum vector
with $\vv={\mm}/{\rho}$ the fluid velocity,
$p$ is the scalar pressure, and
$E$ is the total energy with $e$ the internal energy,
which is
a given function of $(\tau,p)$ or $(\rho,p)$ defined through
thermodynamical relations.
The notation $\aaa\otimes\bb$ denotes the tensor product of
the vectors $\aaa$ and $\bb$.
The other two thermodynamic variables are the temperature $\theta$ and
the entropy $S$. If $(\rho, S)$ are chosen as independent variables,
then the constitutive relations can be written as
$(e,p,\theta)=(e(\rho,S), p(\rho,S),\theta(\rho,S))$
governed by the Gibbs relation:
\begin{equation}
\theta dS=de +pd\tau=de-\frac{p}{\rho^2}d\rho.
\label{1.4}
\end{equation}

For a polytropic gas,
\begin{equation}
p=R\rho\theta=\kappa\rho^\gamma e^{\frac{S}{c_v}},
\quad e=c_v\theta=\frac{\kappa}{\gamma-1}\rho^{\gamma-1}e^{\frac{S}{c_v}}, \qquad \gamma=1+\frac{R}{c_v},
\label{1.5}
\end{equation}
where $R>0$ may be taken as the universal gas constant
divided by the effective molecular weight of the particular gas,
$c_v>0$ is the specific heat at constant volume,
$\gamma>1$ is the adiabatic exponent, and $\kappa>0$ can be
any positive constant by scaling.
The sonic speed is
$c:=\sqrt{p_\rho(\rho,S)}$.
For polytropic gas, $c=\sqrt{\gamma p/\rho}$.

As indicated in \S 2.4 below, no matter how smooth the
initial function is for the Cauchy problem,
the solution of (\ref{Euler1})
generically develops singularities in a finite time.
Then system (\ref{Euler1}) is complemented
by the Clausius-Duhem inequality:
\begin{equation}
\del_t(\rho S)+\nabla_\xx\cdot (\mm S)\ge 0
\label{Entropy1}
\end{equation}
in the sense of distributions
in order to single out physical discontinuous solutions,
so-called {\it entropy solutions}.

\smallskip
When a flow is isentropic (i.e., the entropy $S$ is a uniform constant
$S_0$ in the flow), the Euler equations take
the following simpler form:
\begin{equation}
\left\{\begin{array}{l}
\del_t\rho +\nabla_\xx\cdot \mm=0,\\
\del_t\mm+\nabla_\xx\cdot\left(\frac{\mm\otimes\mm}{\rho}\right)
+\nabla_\xx p=0,
\end{array}\right.
\label{Euler2}
\end{equation}
where the pressure is regarded as a function of density,
$p=p(\rho,S_0)$, with constant $S_0$.
For a polytropic gas,
\begin{equation}
p(\rho)=\kappa\rho^\gamma, \qquad \gamma> 1,
\label{1.7}
\end{equation}
where $\kappa>0$ can be any positive constant under scaling.
When the entropy is initially a uniform constant,
a smooth solution of (\ref{Euler1}) satisfies the equations in (\ref{Euler2}).
Furthermore, it should be observed that the solutions of
system (\ref{Euler2}) are also close to the solutions of
system (\ref{Euler1}) even after shocks form, since the entropy
increases across a shock to third-order in wave strength for solutions
of (\ref{Euler1}),
while in (\ref{Euler2}) the entropy $S$ is constant.
Moreover, system (\ref{Euler2}) is an excellent model for
the isothermal fluid flow with $\gamma=1$ and for the shallow water flow
with $\gamma=2$.

An important example for (\ref{cons-2}) is the potential flow equation:
\begin{equation}\label{potential-1}
\del_t \rho(\Phi_t, \nabla_\xx\Phi)
+\nabla_\xx\cdot \big(\rho (\Phi_t, \nabla_\xx\Phi)\nabla_\xx\Phi\big)
=0,
\end{equation}
where
$\rho(\Phi_t, \nabla_\xx\Phi)
=h^{-1}\big(K-(\del_t\Phi +\frac{1}{2}|\nabla_\xx\Phi|^2)\big),
$
$h'(\rho)=p'(\rho)/\rho$,
for the pressure function $p(\rho)$, which typically takes
form \eqref{1.7} with $\gamma\ge 1$.

The importance of the potential flow equation \eqref{potential-1}
in the time-dependent
Euler flows \eqref{Euler1}
with weak discontinuities
was observed by Hadamard \cite{Hadamard}.
Also see Bers \cite{Bers-book}, Cole-Cook \cite{CoCo}, Courant-Friedrichs \cite{CFr},
Majda-Thomas \cite{MTh}, and Morwawetz \cite{Mora2}.

\section{Basic Features and Phenomena}
We first reveal some basic features and
phenomena of M-D hyperbolic conservation laws
with state domain $\D$, i.e., $\uu\in\D$.

\subsection{Convex Entropy and Symmetrization}

A function $\eta: \D\to \R$ is called an entropy of
system (\ref{cons}) if there exists a vector function
$\qq: \D\to \R^d, \qq=(\qq_1,\dots, \qq_d),$ satisfying
\begin{equation}\label{e-e-f}
\nabla \qq_i(\uu)=\nabla\eta(\uu)\nabla\ff_i(\uu), \qquad i=1,\dots,d.
\end{equation}
Then $\qq$ is called the corresponding entropy flux, and
$(\eta, \qq)$ is simply called an entropy pair.
An entropy $\eta(\uu)$ is called a convex entropy in $\D$
if
$
\nabla^2\eta(\uu)\ge 0$  for any $\uu\in \D$
and a strictly convex entropy in $\D$ if
$
\nabla^2\eta(\uu)\ge c_0 I
$
with a constant $c_0>0$ uniform for $\uu\in \D_1$ for any
$\D_1\subset \bar{\D}_1\Subset \D$, where $I$ is the $m\times m$
identity matrix.
Then the correspondence of the Clausius-Duhem inequality (\ref{Entropy1})
in the context of
hyperbolic conservation laws is the Lax entropy inequality:
\bea
\del_t\eta(\uu)+\nabla_\xx\cdot\qq(\uu)\le 0
\label{Lax}
\eea
in the sense of distributions for any $C^2$ convex entropy
pair $(\eta,\qq)$.

The following important observation
is due to Friedrich-Lax {\rm \cite{FL}},
Godunov {\rm \cite{Go1}}, and Boillat {\rm \cite {Boillat1}}
(also see Ruggeri-Strumia {\rm \cite{Ruggeri-S}}).

\medskip
\begin{theorem} \label{LF-1}
A system in \eqref{cons} endowed with a strictly convex
entropy in $\D$
must be symmetrizable and hence hyperbolic in $\D$.
\end{theorem}

\smallskip
This theorem is particularly useful for determining
whether a large physical system is symmetrizable and hence hyperbolic,
since most of physical systems from continuum
physics are endowed with a strictly convex entropy.
In particular, for system \eqref{Euler1},
\begin{equation}\label{physicalentropy:1}
(\eta_*,\qq_*)=(-\rho S, -\mm S)
\end{equation}
is a strictly convex entropy pair
when $\rho>0$ and $p>0$; while, for system \eqref{Euler2},
the mechanical energy and energy flux pair:
\begin{equation}\label{phyiscalentropy:2}
(\eta_*,\qq_*)=(\rho(\frac{1}{2}\frac{|\mm|^2}{\rho^2}+e),\,
\mm(
\frac{1}{2}\frac{|\mm|^2}{\rho^2}+e(\rho)+\frac{p(\rho)}{\rho}))
\end{equation}
is a strictly convex entropy pair
when $\rho>0$ for polytropic gases.
For a hyperbolic system of conservation laws
without a strictly convex entropy,
it is possible to enlarge the system so that the enlarged system
is endowed with a globally defined, strictly convex entropy.
See Brenier {\rm \cite{Bren}}, Dafermos {\rm \cite{Dafermos-book}},
Demoulini-Stuart-Tzavaras {\rm \cite{DST}},  Qin {\rm \cite{Qin}},
and Serre {\rm \cite{Serre-04}}.

\medskip
There are several direct, important applications of
Theorem \ref{LF-1} based on the symmetry property of
system (\ref{cons}) endowed with a strictly convex entropy.
We list three of them below.

\smallskip
\noindent
{\bf Local Existence of Classical Solutions}.
Consider the Cauchy problem for a general hyperbolic
system (\ref{cons}) with a strictly convex entropy:
\begin{equation} \label{Cauchy}
\uu|_{t=0}=\uu_0.
\end{equation}

\begin{theorem} \label{LS-2}
Assume that $\u_0: \R^d\to \D$ is in $H^s\cap L^{\infty}$ with
$s>d/2+1$.
Then, for the Cauchy problem \eqref{cons} and \eqref{Cauchy},
there exists a finite time
$T=T(\|\u_0\|_s,\|\u_0\|_{L^{\infty}})\in (0,\infty)$
such that there exists a unique bounded classical solution
$\u\in C^1([0,T]\times \R^d)$ with
$\u(t,\xx)\in \D$  for $(t,\xx)\in [0,T]\times \R^d$ and
$
\u\in C([0,T];H^s)\cap C^1([0,T];H^{s-1}).
$
\end{theorem}

\smallskip
Kato \cite{Kato-75} first
formulated and applied a basic idea in the semigroup theory
to yield the local existence of smooth solutions to (\ref{cons}).
Majda in \cite{Ma} provided a proof which relies solely on
the elementary linear existence theory for symmetric hyperbolic
systems with smooth coefficients via a classical iteration
scheme (cf.  \cite{Cou-Hil})
by using the symmetry of system (\ref{cons}).
Moreover, a sharp continuation principle was also provided:
For $\u_0\in H^s$ with $s>d/2+1$,
the interval $[0,T)$ with $T<\infty$ is the maximal
interval of the classical $H^s$ existence for (\ref{cons}) if and only if
either $\u(t,\xx)$ escapes every compact subset $K\Subset\D$ as $t\to T$,
or $\|(\u_t,D\u)(t,\cdot)\|_{L^{\infty}}\to\infty$
as $t\to T$.
The first catastrophe is associated with a blow-up phenomenon
such as focusing and concentration, and the second in this principle
is associated with the formation of shocks
and vorticity waves, among others, in the smooth solutions.
Also see Makino-Ukai-Kawashima \cite{MUK} and Chemin \cite{Chemin-E}
for the local existence
of classical solutions of the Cauchy problem
for the M-D Euler equations.

\smallskip
\noindent
{\bf 2.1.2. Stability of Lipschitz Solutions, Rarefaction Waves, and
Vacuum States in the Class of
Entropy Solutions in $L^\infty$}. Assume that system \eqref{cons} is
endowed with a strictly convex
entropy on compact subsets of $\D$.

\smallskip
\begin{theorem}[Dafermos \cite{Da9,Dafermos-book}]\label{dafermos-1}
Suppose that $\ww$ is a Lipschitz solution of \eqref{cons}
on $[0,T)$, taking values in a convex compact subset $K$ of $\D$,
with initial data $\ww_0$.
Let $\uu$ be any entropy solution of \eqref{cons} on $[0,T)$, taking
values in $K$, with initial data $\uu_0$.
Then
$$
\int_{|\xx|<R}
|\uu(t,\xx)-\ww(t,\xx)|^2d\xx
\le C(T)\int_{|\xx|<R+Lt}
|\uu_0(\xx)-\ww_0(\xx)|^2d\xx
$$
holds for any $R>0$ and $t\in [0,T)$, with $L>0$ depending solely
on $K$ and the Lipschitz constant of $\ww$.
\end{theorem}

\smallskip
The results can be extended to M-D hyperbolic systems of conservation laws
with partially convex entropies and involutions;
see Dafermos {\rm \cite{Dafermos-book}} (also see {\rm \cite{Boi}}).
Some further ideas have been developed \cite{Chen-Chen} to show the
stability of planar rarefaction waves and vacuum states in the
class of entropy
solutions in $L^\infty$ for the M-D Euler equations
by using the theory of divergence-measure fields (see \S 7).

\smallskip
\begin{theorem}[Chen-Chen \cite{Chen-Chen}]\label{Chen-Chen}
Let $\n\in \bS^{d-1}$.
Let
$
{\bf R}(t,\xx)=(\hat{\rho},\hat{\mm})(\frac{\xx\cdot\n}{t})
$
be a planar solution,
consisting of planar rarefaction waves and possible vacuum states,
of the Riemann problem:
$$
{\bf R}|_{t=0}=
(\rho_\pm, \mm_\pm), \qquad \pm \xx\cdot\bnu \ge 0,
$$
with constant states $(\rho_\pm, \mm_\pm)$.
Suppose that $\uu(t,\xx)=(\rho,\mm)(t,\xx)$ is an entropy solution
in $L^\infty$ of
\eqref{Euler2} that may contain vacuum.
Then, for any $R>0$ and $t\in [0,\infty)$,
$$
\int_{|\xx|<R}
\alpha(\uu, {\bf R})(t,\xx)\, d\xx
\le \int_{|\xx|<R+Lt}
 \alpha(\uu, {\bf R})(0,\xx)\, d\xx,
$$
where $L>0$ depends solely on the bounds of
the solutions $\uu$ and ${\bf R}$, and
$
\alpha(\uu, {\bf R})
=(\uu- {\bf R})^\top
(\int_0^1\nabla^2\eta_*({\bf R}+\tau(\uu- {\bf R}))d\tau)
(\uu-{\bf R})
$
with the mechanical energy $\eta_*(\uu)$
in \eqref{phyiscalentropy:2}.
\end{theorem}

\smallskip
A similar theorem to Theorem {\rm \ref{Chen-Chen}} was also established for the adiabatic
Euler equations \eqref{Euler1} with appropriate chosen
entropy function in Chen-Chen {\rm \cite{Chen-Chen}};
also cf. Chen {\rm \cite{Ch6}} and Chen-Frid-Li {\rm \cite{CFL}}.
Also see Perthame \cite{Perth} for the time-decay of the internal energy and the
density for entropy solutions to \eqref{Euler1}
for polytropic gases (also cf. \cite{Chemin-E}).

\smallskip
\noindent
{\bf 2.1.3. Local Existence of Shock-Front Solutions}.
Shock-front solutions, the simplest type of discontinuous solutions,
are the most important discontinuous nonlinear progressing wave
solutions of conservation laws \eqref{cons}.
Shock-front solutions are discontinuous piecewise
smooth entropy solutions with the following structure:
\begin{enumerate}
\item[(i)] There exist a $C^2$ space-time
hypersurface $\mathcal{S}(t)$ defined
in $(t,\xx)$ for $0\le t\le T$ with space-time normal
$({\bf n}_t, {\bf n}_{\xx})=({\bf n}_t,{\bf n}_1,\dots,{\bf n}_d)$
and
two $C^1$ vector-valued functions: $\u^\pm(t,\xx)$,
defined on respective domains $\D^\pm$
on either side of
the hypersurface $\bS(t)$, and satisfying
\begin{equation} \label{SF1}
\partial_t\u^{\pm}+\nabla\cdot\ff(\u^{\pm})=0
\qquad \, \mbox{in}\,\, \D^{\pm};
\end{equation}
\item[(ii)] The jump across $\bS(t)$ satisfies
the Rankine-Hugoniot condition:
\begin{equation} \label{SF2}
\left.\big({\bf n}_t(\u^+-\u^-)
+{\bf n}_{\xx}\cdot(\ff(\u^+)-\ff(\u^-))\big)\right|_{\bS}=0.
\end{equation}
\end{enumerate}
Since \eqref{cons} is nonlinear,
the surface $\bS$ is not known in advance and must be determined
as a part of the solution of the problem; thus the equations in
(\ref{SF1})--(\ref{SF2}) describe a M-D
free-boundary value problem for \eqref{cons}.

The initial functions yielding shock-front solutions are defined
as follows.
Let $\bS_0$ be a smooth hypersurface parameterized by $\a$, and
let ${\bf n}(\a)=({\bf n}_1,\dots, {\bf n}_d)(\a)$ be
a unit normal to $\bS_0$.
Define the piecewise smooth initial functions $\u_0^\pm(\xx)$ for respective
domains $\D_0^\pm$ on either side of the hypersurface
$\bS_0$.
It is assumed that the initial jump satisfies
the Rankine-Hugoniot condition, i.e., there is a smooth scalar
function $\sigma(\a)$ so that
\begin{equation} \label{SF4}
-\sigma(\a)\big(\u_0^+(\a)-\u_0^-(\a)\big)
+{\bf n}(\a)\cdot\big(\ff(\u_0^+(\a))-\ff(\u_0^-(\a))\big)=0
\end{equation}
and that $\sigma(\a)$ does not define a characteristic
direction, i.e.,
\begin{equation} \label{SF5}
\sigma(\a)\neq \lambda_i(\u_0^{\pm}),
\qquad \a\in\overline{\bS_0},\quad 1\le i\le m,
\end{equation}
where $\lambda_i$, $i=1,\dots, m$, are the eigenvalues of (\ref{cons}).
It is natural to require that $\bS(0)=\bS_0$.

\smallskip
Consider the 3-D full Euler equations
in \eqref{Euler1} for $\u=(\rho, \mm, \rho E)$, away from vacuum,
with piecewise smooth initial data:
\begin{equation} \label{euler-ise-i}
\u|_{t=0}=
        \u_0^\pm(\xx), \qquad \xx\in \D_0^\pm.
\end{equation}

\smallskip
\begin{theorem}[Majda \cite{Ma-SF1}]\label{LD1}
Assume that $\bS_0$ is a smooth hypersurface in $\R^3$
and that $\u_0^+(\xx)$ belongs to the uniform
local Sobolev
space $H^s_{ul}(\D_0^+)$, while $\u_0^-(\xx)$
belongs to the Sobolev space $H^s(\D_0^-)$, for some fixed $s\ge 10$.
Assume also that there is a function $\sigma(\a)\in H^s(\bS_0)$ so
that \eqref{SF4}--\eqref{SF5} hold, and the compatibility conditions
up to order $s-1$ are satisfied on $\bS_0$ by the initial data,
together with the entropy condition:
\begin{equation} \label{SF6}
\frac{\mm_0^+}{\r}\cdot{\bf n}(\a)+c(\rho_0^+,S_0^+)<\sigma(\a)<
\frac{\mm_0^-}{\r}\cdot{\bf n}(\a)+ c(\rho_0^-,S_0^-),
\end{equation}
and the Majda stability condition:
\begin{eqnarray}
\Big(\frac{p_S(\rho^+, S^+)}{\rho^+\theta^+}+1\Big)
\frac{p(\rho^+,S^+)-p(\rho^-, S^-)}{\rho^+c^2(\rho^+, S^+)}<1.
\label{SF7}
\end{eqnarray}
Then there exists $T>0$ and a $C^2$-hypersurface $\bS(t)$ together
with $C^1$--functions $\u^\pm(t,\xx)$
defined for $t\in[0,T]$
so that
$\u(t,\xx)=\u^\pm(t,\xx), (t,\xx)\in \D^\pm$,
is the discontinuous shock-front solution
of the Cauchy problem \eqref{Euler1} and \eqref{euler-ise-i}
satisfying \eqref{SF1}--\eqref{SF2}.
\end{theorem}

\smallskip
Condition in \eqref{SF7}
is always satisfied for shocks of any strength
for polytropic gas with $\gamma>1$ and for sufficiently weak shocks for
general equation of state.
In Theorem \ref{LD1}, $c=c(\rho,S)$ is the sonic speed.
The Sobolev space
$H^s_{ul}(\D_0^+)$ is defined as follows:
A vector function $\uu$ is in $H_{ul}^s$, provided that there
exists some $r>0$ such that
$$
\max_{{\bf y}\in \R^d}\|\omega_{r,{\bf y}}\uu\|_{H^s}<\infty
\qquad\mbox{with}\,\,\,
\omega_{r,{\bf y}}(\xx)=\omega(\frac{\xx-{\bf y}}{r}),
$$
where $\omega\in C_0^\infty(\R^d)$ is a function so that $\omega(\xx)\ge 0$,
$\omega(\xx)=1$ when $|\xx|\le \frac{1}{2}$, and $\omega(\xx)=0$ when $|\xx|>1$.

\smallskip
The compatibility conditions in Theorem {\rm \ref{LD1}} are defined in
{\rm \cite{Ma-SF1}} and needed in order to avoid the formation
of discontinuities in higher
derivatives along other characteristic surfaces emanating
from $\bS_0$:
Once the main condition in \eqref{SF4} is satisfied, the compatibility
conditions are automatically guaranteed for a wide class of initial
functions.
Further studies on the local existence and stability of shock-front
solutions can be found in Majda {\rm \cite{Ma-SF1,Ma}}.
The existence of shock-front solutions whose lifespan is uniform with
respect to the shock
strength was obtained in M\'{e}tivier {\rm \cite{Met}}.
Also see Blokhin-Trokhinin {\rm \cite{BT}} for further discussions.

The idea of the proof is similar to that for Theorem {\rm \ref{LS-2}}
by using the existence of a strictly convex
entropy and the symmetrization of \eqref{cons},
but the technical details are quite different due to
the unusual features of the problem
in Theorem {\rm \ref{LD1}}.
For more details, see {\rm \cite{Ma-SF1}}.

The  Navier-Stokes regularization of M-D shocks for the Euler equations
has been established in Gu\`{e}s-M\'{e}tiver-Williams-Zumbrun \cite{GMWZ1}.
The local existence of rarefaction wave front-solutions for M-D
hyperbolic systems of conservation laws has also been
established in Alinhac \cite{Al2}.
Also see Benzoni-Gavage--Serre \cite{BGS}.

\subsection{Hyperbolicity}
There are many examples of $m\times m$ hyperbolic
systems of conservation laws for $d=2$
which are strictly hyperbolic;
that is, they have simple characteristics.
However, for $d=3$, the situation is different.
The following result for the case $m\equiv 2(mod\, 4)$ is due to Lax {\rm \cite{La5}}
and for the more general case
$m\equiv \pm 2,\pm 3, \pm 4 (mod\, 8)$
due to Friedland-Robin-Sylvester {\rm \cite{FRS}}.

\medskip
\begin{theorem}\label{lax:82}
Let $\AAA,\BBB$, and $\CCC$ be three $m\times m$  matrices such that
$
\alpha \AAA +\beta \BBB +\gamma \CCC
$
has real eigenvalues for any real $\alpha,\beta$, and $\gamma$.
If $m\equiv \pm 2,\pm 3, \pm 4\, (mod\, 8)$, then there exist
$\alpha_0,\beta_0,$ and $\gamma_0$ with
$\alpha_0^2+\beta_0^2+\gamma_0^2\ne 0$
such that
$\alpha_0 \AAA +\beta_0 \BBB +\gamma_0 \CCC$
is degenerate, that is, there are at least two eigenvalues
which coincide.
\end{theorem}

This implies that, for $d=3$, there are no strictly hyperbolic systems
if $m\equiv \pm 2,\pm 3, \pm 4\, (mod\, 8)$.
Such multiple characteristics influence the propagation
of singularities.

\medskip
Consider the {\it isentropic Euler equations} (\ref{Euler2}).
When $d=2$ and $m=3$, the system is strictly hyperbolic with
three real eigenvalues $\lambda_-<\lambda_0<\lambda_+$:
$$
\lambda_0=n_1 v_1 + n_2 v_2,
\quad
\lambda_\pm =n_1 v_1+ n_2 v_2\pm c(\rho),
\qquad \rho>0.
$$
Strict hyperbolicity
fails at the vacuum states $\rho=0$.
However, when $d=3$ and $m=4$, the system is no longer strictly
hyperbolic even when $\rho>0$ since the eigenvalue
$
\lambda_0= n_1 v_1 + n_2 v_2+ n_3 v_3
$
has double multiplicity, although the other eigenvalues
$
\lambda_\pm =n_1 v_1+ n_2 v_2+ n_3 v_3\pm c(\rho)
$
are simple when $\rho>0$.

\smallskip
Consider the {\it adiabatic  Euler equations} (\ref{Euler1}).
When $d=2$ and $m=4$, the system is nonstrictly hyperbolic,
since the eigenvalue
$
\lambda_0=n_1 v_1 + n_2 v_2
$
has double multiplicity, although
$
\lambda_\pm =n_1 v_1+ n_2 v_2\pm c(\rho,S)
$
are simple when $\rho>0$.
When $d=3$ and $m=5$, the system is again nonstrictly hyperbolic,
since the eigenvalue
$
\lambda_0=n_1 v_1 +n_2 v_2+n_3 v_3
$
has triple multiplicity, although
$
\lambda_\pm =n_1 v_1+n_2 v_2+n_3 v_3
\pm c(\rho,S)
$
are simple when $\rho>0$.

\subsection{Genuine Nonlinearity}
The $j^{th}$-characteristic field of system (\ref{cons})
in $\D$
is called {\it genuinely nonlinear} if,
for each fixed $\bnu\in S^{d-1}$,
the $j^{th}$-eigenvalue $\lambda_j(\uu;\bnu)$
and the corresponding eigenvector $\rr_j(\uu;\bnu)$
determined by
$
(\nabla\ff(\uu)\cdot \bnu)\rr_j(\uu;\bnu)
=\lambda_j(\uu;\bnu)\rr_j(\uu;\bnu)
$
satisfy
\begin{equation}\label{2.4.1}
\nabla_{\uu}\lambda_j(\uu;\bnu)\cdot
\rr_j(\uu;\bnu)\ne 0
\qquad \mbox{for any}\,\, \uu\in \D, \, \bnu\in \bS^{d-1}.
\end{equation}
The $j^{th}$-characteristic field of (\ref{cons})
in  $\D$ is called linearly degenerate if
\begin{equation}\label{2.4.2}
\nabla_{\uu}\lambda_j(\uu;\bnu)\cdot
\rr_j(\uu;\bnu)\equiv 0
\qquad \mbox{for any}\,\, \uu\in \D.
\end{equation}

Then any {\it scalar} quasilinear conservation law in $\R^d, d\ge 2,$
is never genuinely nonlinear in all directions.
This is because, in this case,
$\lambda(u;\bnu)=\ff''(u)\cdot\bnu$ and $r(u;\bnu)=1$, and
$
\lambda'(u;\bnu) r(u;\bnu)=\ff''(u)\cdot\bnu,
$
which is impossible to make this never equal to zero
in all directions.
A M-D version of
genuine nonlinearity for scalar conservation laws is
that the set

\smallskip
\noindent
\quad
{\it $\{u\, :\, \tau +\ff'(u)\cdot \bnu=0\}$
has zero Lebesgue measure for any $(\tau, \bnu)\in \bS^d$,}

\smallskip
\noindent
which is a generalization of \eqref{2.4.1}.
Under this generalized nonlinearity, the following have been
established:
(i) Solution operators are compact as $t>0$
in Lions-Perthame-Tadmor \cite{LPT1} (also see \cite{CF5,RT});
(ii) Decay of periodic solutions (Chen-Frid \cite{CF2}; also see
Engquist-E  \cite{EE});
(iii) Strong initial and boundary traces of entropy solutions
      (Chen-Rascle \cite{CR}, Vasseur \cite{Vas};
      also see Panov \cite{panov});
(iv) $BV$-like structure of $L^\infty$ entropy solutions
(De Lellis-Otto-Westdickenberg \cite{DOW}).
Furthermore, we have

\smallskip
\begin{theorem}[Lax \cite{La5}]\label{Lax84}
Every real, strictly hyperbolic quasilinear
system for $m=2k$, $k\ge 1$ odd, and $d=2$
is linearly degenerate in some direction.
\end{theorem}

\smallskip
Quite often, linear degeneracy results from the loss of strict
hyperbolicity. For example, even in the 1-D case,
if there exists $i\ne j$ such that
$\lambda_i(\uu)=\lambda_j(\uu)$ for all $\uu\in K$,
then Boillat {\rm \cite{Boi}} proved that
the $i^{th}$- and $j^{th}$-characteristic fields are linearly
degenerate.

\medskip
For the {\it isentropic Euler equations} (\ref{Euler2}) with
$d=2$ and $m=3$ for polytropic gases with the eigenvalues:
\begin{eqnarray*}
\lambda_0=n_1 v_1+n_2 v_2,
\qquad
\lambda_\pm=n_1 v_1+ n_2 v_2\pm c(\rho),
\end{eqnarray*}
and the corresponding eigenvectors $\rr_0$ and $\rr_\pm$,
we have
$
\nabla \lambda_0\cdot\rr_0\equiv 0,
$
which is linearly degenerate,
and
$\nabla \lambda_\pm\cdot\rr_\pm\ne 0$ for $\rho\in (0, \infty)$,
which are genuinely nonlinear.
For the {\it adiabatic Euler equations} (\ref{Euler1}) with $d=2$ and $m=4$
for polytropic gases with the eigenvalues:
\begin{eqnarray*}
\lambda_0=n_1 v_1+ n_2 v_2,
\qquad
\lambda_\pm= n_1 v_1+ n_2 v_2\pm c(\rho,S),
\end{eqnarray*}
and the corresponding eigenvectors $\rr_0$ and $\rr_\pm$,
we have
$
\nabla \lambda_0\cdot\rr_0\equiv 0,
$
which is linearly degenerate, and
$\nabla \lambda_\pm\cdot\rr_\pm\ne 0$ for $\rho, S\in (0, \infty)$,
which are genuinely nonlinear.

\subsection{Singularities}
For the 1-D case, singularities
include the formation of shocks, contact discontinuities, and
the development of vacuum states, at least for gas dynamics.
For the M-D case, the situation is much
more complicated: besides shocks, contact discontinuities, and vacuum states,
singularities may include vorticity waves,
focusing waves, concentration waves,
complicated wave interactions, among others.
Consider the Cauchy problem of the Euler equations
in (\ref{Euler1}) in $\R^3$ for polytropic gases with smooth
initial data:
\begin{equation} \label{3d-euler-i}
\u|_{t=0}=\u_0(\xx) \qquad
\mbox{with}\,\, \r_0(\xx)>0,
\end{equation}
satisfying
$(\r_0,\mm_0,S_0)(\xx)=(\bar{\r},0,\bar{S})$ for $|\xx|\ge L$,
where $\bar{\r}>0$, $\bar{S}$, and $L$ are constants.
The equations in (\ref{Euler1}) possess a unique local $C^1$-solution
$\u(t,\xx)$ with $\r(t,\xx)>0$ provided that
the initial function (\ref{3d-euler-i}) is sufficiently regular
(Theorem \ref{LS-2}).
The support of the smooth disturbance
$(\r_0(\xx)-\bar{\r},\, \mm_0(\xx),\, S_0(\xx)-\bar{S})$ propagates with
speed at most
$\bar{c}=\sqrt{p_{\r}(\bar{\r},\bar{S})}$ (the sound speed)
that is,
$(\r,\mm,S)(t,\xx)=(\bar{\r},0,\bar{S})$
if $|\xx|\ge L+\bar{c} t$.
Take $\bar{p}=p(\bar{\r},\bar{S})$. Define
\begin{eqnarray*}
P(t)=\int_{{\R}^3}\big(p(t,\xx)^{1/\gamma}-\bar{p}^{1/\gamma}\big)d\xx,\qquad
F(t)=\int_{{\R}^3}\mm(t,\xx)\cdot\xx \, d\xx,
\end{eqnarray*}
which measure the entropy and the radial component
of momentum.

\medskip
\begin{theorem}[Sideris \cite{S1}] \label{blowup-3d}
Suppose that $(\r,\mm,S)(t,\xx)$ is a $C^1$--solution of \eqref{Euler1}
and \eqref{3d-euler-i} for $0<t<T$ and
\begin{eqnarray}
P(0)\ge 0, \qquad
F(0)> \frac{16}{3}\pi \bar{c}L^4\max_{\xx}\big(\r_0(\xx)\big).
\label{blowup-2}
\end{eqnarray}
Then the lifespan $T$ of the $C^1$--solution is finite.
\end{theorem}

In particular, when $\r_0=\bar{\r}$ and $S_0=\bar{S}$,
$P(0)=0$ and the second condition in \eqref{blowup-2} holds
if the initial velocity
$\v_0(\xx)=\frac{\mm_0(\xx)}{\rho_0(\xx)}$ satisfies
$
\int_{|\xx|<L}\v_0(x) \cdot\xx \, d\xx > \bar{c}\sigma L^4.
$
From this, one finds that the initial velocity must be supersonic
in some region relative to $\bar{c}$.
The formation of a singularity (presumably a shock) is detected
as the disturbance overtakes the wave-front forcing the front to propagate
with supersonic speed.
The formation of singularities occurs even without condition of
largeness such as \eqref{blowup-2}. For more details, see \cite{S1}.

In Christodoulou {\rm \cite{Christodoulou1}},  the relativistic Euler equations
for a perfect fluid with an arbitrary equation of state have been analyzed.
The initial function is imposed on a given
spacelike hyperplane and is constant outside a compact set.
Attention is restricted to the evolution
of the solution within a region limited by two concentric spheres.
Given a smooth solution, the geometry of the boundary of its domain
of definition is studied, that is, the locus where shocks may form.
Furthermore, under certain smallness assumptions on the size of the initial data,
a remarkable and complete picture of the formation of shocks in $\R^3$
has been obtained. In addition, sharp sufficient conditions on the initial
data for the formation of shocks in the evolution have been identified,
and sharp lower and upper bounds for the time
of existence of a smooth solution have been derived.
Also see Christodoulou {\rm \cite{Christodoulou2}}
for the formation of black holes in general relativity.

\subsection{$BV$ Bound}
For 1-D strictly hyperbolic systems,
Glimm's theorem \cite{Glimm}
indicates that, as long as $\|\uu_0\|_{BV}$ is
sufficiently small, the solution $\uu(t,x)$ satisfies
the following stability estimate:
\begin{equation}\label{BV-1}
 \|\uu(t,\cdot)\|_{BV}\le C\, \|\uu_0\|_{BV}.
\end{equation}
Even more strongly,
for two solutions $\uu(t,x)$ and $\ww(t,x)$ obtained by either
the Glimm scheme, wave-front tracking method, or vanishing
viscosity method with small total variation,
$$
\|\u(t,\cdot)-\ww(t,\cdot)\|_{L^1(\R)}
\le C\,\|\ww(0,\cdot)-\v(0,\cdot)\|_{L^1(\R)}.
$$
See Bianchini-Bressan \cite{BBr2} and Bressan \cite{Bre-book};
also see
Dafermos \cite{Dafermos-book},
Holden-Risebro \cite{HRis}, LeFloch \cite{LeF1}, Liu-Yang \cite{LY-2},
and the references cited therein.

The recent great progress for entropy solutions to
1-D hyperbolic conservation laws
based on $BV$ estimates and trace theorems of $BV$ fields
naturally arises the expectation that a similar approach may also be
effective for M-D hyperbolic
systems of conservation laws, that is,
whether entropy solutions satisfy the relatively modest stability
estimate \eqref{BV-1}.
Unfortunately, this is not the case.
Rauch \cite{Rauch86} showed that
the necessary condition for (\ref{BV-1})
to be held is
\begin{equation}\label{bound:4}
\nabla \ff_k(\uu)\nabla \ff_l(\uu)
=\nabla\ff_l(\uu) \nabla\ff_k(\uu)
\qquad \mbox{for all}\,\, k, l=1,2,\dots,d.
\end{equation}
The analysis above suggests that only systems in which the
commutativity
relation (\ref{bound:4}) holds offer some hope for treatment in
the $BV$ framework.
This special case includes the scalar case $m=1$ and the 1-D
case $d=1$.
Beyond that, it contains very few systems of physical interest.
An example is the system with fluxes:
$
\ff_k(\uu)=\phi_k(|\uu|^2)\uu$,
$k=1,2,\dots, d$,
which governs the flow of a fluid in an anisotropic porous
medium.
However, the recent study in
Bressan \cite{Bressan03}
and Ambrosio-DeLellis \cite{Am-De1} shows that, even in this case,
the space $BV$ is not a well-posed space
(also cf. Jenssen \cite{Jen}).
Moreover,
entropy solutions generally do not have even the
relatively modest stability:
$\|\uu(t,\cdot)-\bar{\uu}\|_{L^p}\le C_p\|\uu_0-\bar{\uu}_0\|_{L^p}$,
$p\ne 2$,
based on the linear theory by Brenner \cite{Bre}.

In this regard, it is important to identify
good analytical frameworks for the analysis of entropy solutions of
M-D conservation laws (\ref{cons}),
which are not in $BV$, or even in $L^p$.
A general framework is the space
of divergence-measure fields, formulated recently
in \cite{CF6,ChenTorres,ChenTorresZiemer1},
which is based on
entropy solutions satisfying the Lax entropy inequality
\eqref{Lax}. See \S 7 for more details.

Another important notion of solutions is the notion of the measure-valued entropy solutions,
based on the Young measure representation of a weak convergent sequence
(cf. \cite{Ta1,Ba2}),
proposed by DiPerna \cite{DiPerna-YM}.
An effort has been made to establish the existence of measure-valued
solutions to the M-D Euler equations for compressible fluids;
see Ganbo-Westdickenberg \cite{GWest}.

\subsection{Nonuniqueness}
Another important feature is the nonuniqueness of entropy solutions in general.
In particular,  De Lellis-Sz\'{e}kelyhidi \cite{DeLS1} recently showed
the following
remarkable fact.

\smallskip
\begin{theorem}
Let $d\ge 2$. Then, for any given function $p=p(\rho)$ with
$p'(\rho)>0$ when $\rho>0$, there exist bounded initial functions
$(\rho_0, \vv_0)=(\rho_0, \frac{\mm_0}{\rho_0})$ with $\rho_0(\xx)\ge \delta_0>0$ for which there exist
infinitely many bounded solutions $(\rho, \mm)$ of \eqref{Euler2}
with $\rho\ge \delta$ for some $\delta>0$,
satisfying the energy identity in the sense of distributions:
\begin{equation}\label{energy-identity}
\partial_t (\rho E)
+\nabla_\xx\cdot \big(\mm(E+p/\rho)\big)=0.
\end{equation}
\end{theorem}
In fact, the same result also holds for the full
Euler system, since the solutions constructed satisfy the energy
equality so that there is no energy production at all.
The main point for the result is that the solutions constructed
contain only vortex sheets and vorticity waves which keep
the energy identity \eqref{energy-identity} even for weak solutions
in the sense of distributions, while the vortex sheets do not
appear in the 1-D case.
The construction of the infinitely many solutions is based on
a variant of the Baire category method for differential inclusions.
Therefore, for the uniqueness issue, we have to narrow down further
the class of entropy solutions to single out physical relevant solutions,
at least for the Euler equations.

\section{Multidimensional Systems and Models}
M-D problems are extremely
rich and complicated.
Some great developments have been made
in the recent decades through strong and close interdisciplinary
interactions and diverse approaches, including

\begin{enumerate}
\item[(i)] Experimental data;

\item[(ii)] Large and small scale computing by a search
    for effective numerical methods;

\item[(iii)] Asymptotic and qualitative modeling;

\item[(iv)] Rigorous proofs for prototype problems
    and an understanding of the solutions.
\end{enumerate}
In some sense, the developments
made by using approach (iv) are still behind those by using
the other approaches (i)--(iii);
however, most scientific problems
are considered to be solved satisfactorily only after
approach (iv) is achieved.
In this section, together with Sections 4--7, we present some
samples of M-D systems/models and related important problems
with emphasis on those
prototypes that have been solved or may be expected to
be solved rigorously at least for some cases.

Since the M-D problems are so complicated in general,
a natural strategy to attack these problems as a first
step is to
study (i) simpler nonlinear systems with strong
   physical motivations and (ii) special, concrete nonlinear
   physical problems.
Meanwhile, extend the results and ideas from the first step to
(i) study the Euler equations in gas dynamics and elasticity;
(ii) study more general problems;
(iii) study nonlinear systems that the Euler equations are
the main subsystem or describe the dynamics of macroscopic variables
such as
Navier-Stokes equations,
MHD equations, combustion equations,
Euler-Poisson equations,
kinetic equations especially including the Boltzmann equation,
among others.

Now we first focus on some samples of M-D
systems and models.

\subsection{Unsteady Transonic Small Disturbance Equation}

A simple model
in transonic aerodynamics is the UTSD equation or so-called the 2-D inviscid
Burgers equation
(see Cole-Cook \cite{CoCo}):
\begin{equation}\label{Euler11}
\left\{\begin{array}{l}
\del_t u +\del_x (\frac{1}{2}u^2)+\del_y v=0,\\
\del_y u-\del_x v=0,
\end{array} \right.
\end{equation}
or in the form of Zabolotskaya-Khokhlov equation \cite{ZK}:
\begin{equation}\label{Euler12}
\del_t(\del_tu +u\del_x u) +\del_{yy} u=0.
\end{equation}

The equations in (\ref{Euler11}) describe the potential flow field near the
reflection point in weak
shock reflection, which determines the leading order approximation of
geometric optical expansions; and it can also be used to formulate
asymptotic equations for the transition from regular to Mach
reflection for weak shocks. See Morawetz \cite{Mora2},
Hunter \cite{hunter1}, and the references
cited therein.
Equation (\ref{Euler12}) also arises in many different situations;
see \cite{hunter1,Timm,ZK}.

\subsection{Pressure-Gradient Equations and Nonlinear Wave Equations}
The inviscid fluid motions are driven mainly by the pressure  gradient
and the fluid convection (i.e., transport).
As for modeling, it is natural to study
first the effect of the two driving factors separately.

Separating the pressure gradient from the Euler equations \eqref{Euler1},
we first have the pressure-gradient system:
\begin{equation}\label{Euler13}
\del_t \rho=0,\quad
\del_t(\rho \vv) +\nabla_\xx p =0,\quad
\del_t(\rho e) + p\nabla_\xx\cdot\vv=0.
\end{equation}
We may choose $\rho=1$.
Setting $p=(\gamma-1)P$ and $t=\frac{s}{\gamma-1}$,
and eliminating the velocity $\vv$,
we obtain the following nonlinear wave equation for $P$:
\begin{equation}\label{Euler14}
\del_{ss}(\ln P) -\Delta_\xx P=0.
\end{equation}
Although system (\ref{Euler13}) is obtained from
the splitting idea,
it is a good approximation
to the full Euler equations, especially when the velocity
$\vv$ is small and the adiabatic exponent $\gamma>1$
is large. See Zheng \cite{Zhe2}.

Another related model is the following nonlinear wave equation proposed by
Canic-Keyfitz-Kim \cite{CKK06}:
\begin{equation}\label{nlw}
\del_{tt}\r -\Delta_\xx p(\r)=0,
\end{equation}
where $p=p(\r)$ is the pressure-density relation for fluids.
Equation \eqref{nlw} is obtained from \eqref{Euler2} by neglecting
the inertial terms, i.e., the quadratic terms in $\mm$. This yields
the following system:
\begin{equation}\label{nlw-2}
\del_t\r +\nabla_\xx\cdot \mm=0,\quad
\del_t\mm +\nabla_\xx p(\rho)=0,
\end{equation}
which leads to \eqref{nlw} by eliminating $\mm$ in the system.

\subsection{Pressureless Euler Equations}
$\,\,$ With the pressure-gradient equations (\ref{Euler13}),
the convection (i.e., transport)
part of fluid flow forms the pressureless Euler equations:
\begin{equation}\label{Euler15}
\left\{\begin{array}{ll}
\del_t \rho +\nabla_\xx \cdot(\rho \vv)=0,\\
\del_t(\rho \vv) +\nabla_\xx\cdot (\rho \vv\otimes \vv)=0,\\
\del_t (\rho E) +\nabla_\xx\cdot (\rho E\vv)=0.
\end{array}\right.
\end{equation}
This system also models the motion of free particles which stick
under collision;
see Brenier-Grenier \cite{BG}, E-Rykov-Sinai \cite{ERS},
and Zeldovich \cite{Zel}.
In general, solutions of (\ref{Euler15}) become measure solutions.

System (\ref{Euler15}) has been analyzed extensively; cf.
\cite{B,BG,CLiu,ERS,Gre,HW,Liq,
SZ,WHD} and the
references cited therein.
In particular, the existence of  measure solutions
of the Riemann problem was first presented in Bouchut \cite{B}
for the 1-D case.
It has been shown that
$\delta$-shocks and vacuum states do occur in the Riemann
solutions even in the 1-D case.
Since the two eigenvalues of the transport equations coincide,
the occurrence as $t>0$
can be regarded as a result of resonance between the two
characteristic fields.
Such phenomena can also be regarded as
the phenomena of concentration and cavitation
in solutions as the pressure vanishes.
It has been rigorously shown in Chen-Liu \cite{CLiu} for $\gamma>1$
and Li \cite{Liq} for $\gamma=1$
that, as the pressure vanishes, any two-shock Riemann
solution to the Euler equations tends to
a $\delta$-shock solution to (\ref{Euler15})
and the intermediate densities between the two shocks tend to a weighted
$\delta$-measure that forms the $\delta$-shock.
By contrast, any two-rarefaction-wave Riemann solution of
the Euler equations
has been shown in \cite{CLiu} to tend to a two-contact-discontinuity
solution to (\ref{Euler15}),
whose intermediate state between the two contact discontinuities
is a vacuum state,
even when the initial density stays away from the vacuum.

\subsection{Incompressible Euler Equations}
The incompressible
Euler equations take the form:
\begin{equation}\label{Euler10}
\left\{\begin{array}{ll}
\del_t \vv +\nabla\cdot(\vv\otimes \vv) +\nabla p=0, \\
\nabla\cdot \vv =0
\end{array}\right.
\end{equation}
for the density-independent case,
and \begin{equation}\label{incom:2}
\left\{\begin{array}{ll}
\del_t\rho +\nabla\cdot(\rho \vv)=0,\\
\del_t (\rho \vv) +\nabla\cdot(\rho \vv\otimes \vv) +\nabla p=0, \\
\nabla\cdot \vv =0
\end{array} \right.
\end{equation}
for the density-dependent case, where $p$ should be regarded as an unknown.
These systems can be obtained by formal asymptotics for
low Mach number expansions from
the Euler equations \eqref{Euler1}.
For more details, see
Chorin \cite{Chorin},
Constantin \cite{Cons}, Hoff \cite{Hoff},
Lions \cite{Lio}, Lions-Masmoudi \cite{LM},
Majda-Bertozzi \cite{MB},
and the references cited therein.

System \eqref{Euler10} or \eqref{incom:2}
is an excellent model to capture M-D
vorticity waves by ignoring the shocks in fluid flow,
while system (\ref{potential-1})
is an excellent model to capture M-D shocks by
ignoring the vorticity waves.

\subsection{Euler Equations in Nonlinear Elastodynamics}
\label{elastodynamics}

The equations of nonlinear elastodynamics provide another
excellent example of the rich special structure one encounters
when dealing with hyperbolic
conservation laws.
In $\R^3$, the state vector is $(\vv,\FF)$, where
$\vv\in \R^3$ is the velocity vector and $\FF$ is
the $3\times 3$ matrix-valued deformation gradient
constrained by the requirement that $\det \FF>0$.
The system of conservation laws,
which express the integrability conditions between $\vv$ and $\FF$
and the balance of linear momentum, reads
\begin{equation}\label{NE}
\left\{\begin{array}{ll}
\partial_t F_{i\alpha}-\partial_{x_\alpha}v_i =0,
\qquad \qquad\quad \, i,\alpha=1,2,3,\\
\partial_t v_j -\sum_{\beta=1}^{3}\partial_{x_\beta}S_{j\beta}(F)=0,
\quad j=1,2,3.
\end{array}\right.
\end{equation}
The symbol $\BS$ stands for the {\em Piola-Kirchhoff stress tensor},
which is determined by the (scalar-valued)
{\em strain energy function} $\sigma (\FF)$:
$
S_{j\beta}(\FF)=\frac{\partial\sigma(\FF)}{\partial F_{j\beta}}.
$
System (\ref{NE}) is hyperbolic if and only if, for any
vectors ${\bf \xi}, {\bf n}\in \BS^3$,
\begin{equation} \label{NE2}
\sum_{1\le i,j,\alpha,\beta\le 3}
\;
\frac{\partial^2 \sigma(\FF)}{\partial F_{i\alpha}\partial F_{j\beta}}\;
\xi_i \xi_j n_{\alpha}n_{\beta}>0.
\end{equation}
System (\ref{NE}) is endowed with an entropy pair:
$$
\eta=\sigma(\FF) +\frac12|\vv|^2 ,
\qquad q_{\alpha}=-\sum_{1\le j\le 3} v_j S_{j\alpha}(\FF).
$$
However, the laws of physics do not allow $\sigma(\FF)$, and thereby $\eta$,
to be convex functions.  Indeed, the convexity  of $\sigma$ would violate
the principle of {\em material frame indifference}
$\sigma(\OO\FF)=\sigma (\FF)$ for all $\OO\in SO (3)$
and would also be incompatible with the natural requirement that
$\sigma(\FF)\to\infty$ as $\det \FF\downarrow 0$
or $\det \FF\uparrow \infty$ (see Dafermos \cite{Da-86}).

The failure of $\sigma$ to be convex is also the main source
of complication in elastostatics, where one is seeking to determine
equilibrium configurations
of the body by minimizing the total strain energy $\int\sigma(\FF) d\xx$.
The following  alternative conditions, weaker than convexity
and physically reasonable,
are relevant in that context:

\smallskip
\begin{enumerate}
\item[(i)] {\em Polyconvexity} in the sense of Ball \cite{Ba1}:
$
\sigma(\FF)=g(\FF,\FF^{\ast},\det \FF),
$
where $\FF^{\ast}=(\det \FF)\FF^{-1}$ is the adjugate
of $\FF$ (the matrix of cofactors of $\FF$)
and $g(\FF,\GG,w)$ is a convex function of $19$ variables;

\item[(ii)] {\em Quasiconvexity} in the sense of Morrey
\cite{Morry52};

\item[(iii)] {\em Rank-one convexity}, expressed by (\ref{NE2}).
\end{enumerate}

\noindent
It is known that convexity $\Rightarrow$
polyconvexity $\Rightarrow$
quasiconvexity $\Rightarrow$
rank-one convexity,
however, none of the converse statements is generally valid.

It is important to investigate the relevance of the above
conditions in elastodynamics.
A first start was made in Dafermos \cite{Da-86}
where it was shown that rank-one convexity suffices for the
local existence of classical solutions, quasiconvexity yields
the uniqueness of classical solutions
in the class of entropy solutions in $L^\infty$,
and polyconvexity renders the system symmetrizable
(also see \cite{Qin}).
To achieve this for polyconvexity, one of the main
ideas is to enlarge system (\ref{NE})
into a large, albeit
equivalent,
system for the new state vector
$(\vv,\FF,\FF^*, w)$ with $w={\rm det}\, \FF$:
\begin{eqnarray}
&&\del_t w=\sum_{1\le\alpha,i\le 3}
              \del_{x_\alpha}(F^*_{\alpha i} v_i),
              \label{Na-1}\\
&&\del_t F^*_{\gamma k}
  =\sum_{1\le\alpha,\beta, i,j\le 3}
   \del_{x_\alpha}(\epsilon_{\alpha\beta\gamma}
        \epsilon_{ijk}F_{j\beta}v_i),
   \quad \gamma,k=1,2,3, \label{Na-2}
\end{eqnarray}
where $\epsilon_{\alpha\beta\gamma}$ and $\epsilon_{ijk}$
denote the standard permutation symbols.
Then the enlarged system (\ref{NE}) and (\ref{Na-1})--(\ref{Na-2})
with $21$ equations is endowed a uniformly
convex entropy
$
\eta=\sigma(\FF,\FF^*,w)+\frac{1}{2}|\vv|^2
$
so that the local existence of classical solutions
and the stability of Lipschitz solutions may be inferred directly
from Theorem 2.3.
See Dafermos \cite{Dafermos-book},
Demoulini-Stuart-Tzavaras \cite{DST}, and Qin \cite{Qin}
for more details.

\subsection{Born-Infeld System in Electromagnetism}

The Born-Infeld system is a nonlinear version of the Maxwell
equations:
\begin{equation}\label{BI}
\partial_{t}B+{\rm curl}\big(\frac{\partial W}{\partial D}\big)=0,\quad
\partial_{t}D-{\rm curl}\big(\frac{\partial W}{\partial B}\big)=0,
\end{equation}
where $W:\R^{3}\times\R^{3}\rightarrow\R$ is the given energy density.
The Born--Infeld model corresponds to the special case
$$
W_{BI}(B,D)=\sqrt{1+|B|^{2}+|D|^{2}+|P|^{2}}.
$$
When $W$ is strongly convex (i.e., $D^2W>0$), system (\ref{BI})
is endowed with a strictly convex entropy.
However, $W_{\rm BI}$ is not convex for a large enough field.
As in \S \ref{elastodynamics},
the Born-Infeld model
is enlarged from $6$ to $10$ equations in Brenier \cite{Bren},
by an adjunction of the conservation laws satisfied
by $P:= B\times D$ and $W$ so that
the augmented system turns out to be a set of conservation laws
in the unknowns
$
(h,B,D,P)\in\R\times\R^{3}\times\R^{3}\times\R^{3},
$
which is in the physical region:
$$
\{(h,B,D,P)\,:\,P=D\times B, h=\sqrt{1+|B|^{2}+|D|^{2}+|P|^{2}}>0\}.
$$
Then the enlarged system is
\begin{equation}\label{BI}
\left\{\begin{array}{ll}
\partial_{t}h+{\rm div}P  =  0,  \\
    \smallskip
\partial_{t}B+{\rm curl}\left(\frac{P\times B+D}{h}\right)= 0,  \\
 \smallskip
\partial_{t}D+ {\rm curl}\left(\frac{P\times D-B}{h}\right)= 0,  \\
 \smallskip
\partial_{t}P+{\rm Div}\left(\frac{P\otimes P-B\otimes B-D\otimes D-I}{h}
\right)= 0,
\end{array}\right.
\end{equation}
which is endowed with a strongly convex entropy,
where $I$ is the $3\times 3$ identity matrix.
Also see Serre \cite{Serre-04}
for another enlarged system consisting of 9 scalar evolution
equations in 9 unknowns $(B,D,P)$, where $P$ stands for the
relaxation of the expression $D\times B$.

\subsection{Lax Systems}
Let $f(\uu)$ be an analytic function of a single complex variable
$\uu=u+ v i$.
We impose on the complex-valued function
$\uu=\uu(t,z), z=x+y i,$ and the real variable $t$
the following nonlinear PDE:
\begin{equation}\label{4.6.3}
\del_t \bar{\uu}+\del_z f(\uu)=0,
\end{equation}
where the bar denotes the complex conjugate
and $\del_z=\frac{1}{2}(\del_x-i\del_y)$.
We may express this
equation in terms of the real and imaginary parts of $\uu$
and
$
\frac{1}{2}f(\uu)=a(u,v)+ b(u,v) i.
$
Then (\ref{4.6.3}) gives
\begin{equation}\label{4.6.3a}
\left\{\begin{array}{ll}
\del_tu+\del_x a(u,v) +\del_y b(u,v)=0,\\
\smallskip
\del_tv -\del_x b(u,v) +\del_y a(u,v)=0.
\end{array}\right.
\end{equation}
In particular, when $f(\uu)=\uu^2=u^2+v^2+2uv i$,
system (\ref{4.6.3}) is called the complex Burger equation,
which becomes
\begin{equation}\label{4.6.3b}
\left\{\begin{array}{ll}
\del_tu+ \frac{1}{2}\del_x (u^2+v^2) +\del_y (uv)=0,\\
\smallskip
\del_tv -\del_x (uv) +\frac{1}{2}\del_y (u^2+v^2)=0.
\end{array}\right.
\end{equation}
System (\ref{4.6.3a}) is a symmetric hyperbolic system
of conservation laws with a strictly convex entropy
$
\eta(u,v)= u^2+v^2;
$
see Lax \cite{La6} for more details.
For the 1-D case, this system
is an archetype of hyperbolic conservation
laws with umbilic degeneracy, which
has been
analyzed in Chen-Kan \cite{CK}, Schaeffer-Shearer \cite{SS1},
and the references cited therein.

\subsection{Gauss-Codazzi System for Isometric Embedding}
A fundamental problem in differential geometry is to characterize
intrinsic metrics on a 2-D Riemannian manifold
${\mathcal M}^2$ which can be realized as isometric immersions into
$\R^3$ (cf. Yau \cite{Yau00}; also see \cite{HanHong}).
For this, it suffices to solve the Gauss-Codazzi system,
which can be written as
(cf. \cite{CSW2,HanHong})
\begin{equation} \label{g1}
\left\{\begin{split}
\d_x{M}-\d_y{L}&=\Gamma^{(2)}_{22}L-2\Gamma^{(2)}_{12}M+\Gamma^{(2)}_{11}N, \\
\d_x{N}-\d_y{M}&=-\Gamma^{(1)}_{22}L+2\Gamma^{(1)}_{12}M-\Gamma^{(1)}_{11}N,
\end{split}\right.
\end{equation}
and
\begin{equation}\label{g2}
LN-M^2=K.
\end{equation}
Here
$K(x,y)$ is the given Gauss
curvature, and
$
\Gamma_{ij}^{(k)}
$
is the given Christoffel symbol.
We now follow Chen-Slemrod-Wang \cite{CSW2}
to present the fluid dynamical formulation of system \eqref{g1}--\eqref{g2}.
Set
$
L=\r v^2+p, M=-\r uv, N=\r u^2+p,
$
and  set $q^2=u^2+v^2$ as usual. Then the Codazzi equations \eqref{g1}
become the familiar balance laws of momentum:
\begin{equation} \label{g3}
\left\{\begin{split}
&\d_x(\r uv)+\d_y(\r v^2+p)
 =-\Gamma^{(2)}_{22}(\r v^2+p) -2\Gamma^{(2)}_{12}\r uv-\Gamma^{(2)}_{11}(\r u^2+p), \\
&\d_x(\r u^2+p)+\d_y(\r uv)
 =-\Gamma^{(1)}_{22}(\r v^2+p)-2\Gamma^{(1)}_{12}\r uv-\Gamma^{(1)}_{11}(\r u^2+p),
\end{split}\right.
\end{equation}
and the Gauss equation \eqref{g2} becomes
$
\r p q^2+p^2=K.
$
We choose $p$ to be the
Chaplygin gas-type to allow the Gauss curvature $K$ to change sign:
$p=-\frac{1}{\r}$.
Then we
obtain the ``Bernoulli" relation:
\begin{equation}\label{g6}
\r=\frac{1}{\sqrt{q^2+K}},
\qquad\mbox{or}\qquad
p=-\sqrt{q^2+K}.
\end{equation}
In general, system \eqref{g3}--\eqref{g6} for
unknown $(u,v)$ is of mixed
hyperbolic-elliptic type
determined by the sign of Gauss curvature (hyperbolic when $K<0$
and elliptic when $K>0$).

For the Gauss-Codazzi-Ricci equations for isometric embedding
of higher dimensional Riemannian manifolds,
see Chen-Slemrod-Wang \cite{CSW3}.

\section{Multidimensional Steady Supersonic Problems}

M-D steady problems for the Euler equations are
fundamental in fluid mechanics. In particular, understanding
of these problems helps us to understand the asymptotic
behavior of evolution solutions for large time,
especially global attractors.
One of the excellent sources of steady problems is
Courant-Friedrichs's book \cite{CFr}.
In this section we first discuss some of recent developments
in the analysis of M-D steady supersonic problems.
The M-D steady Euler flows are governed by
\begin{equation}\label{eq1.1}
\left\{%
\begin{array}{ll}
\nabla_\xx\cdot (\rho \vv)=0,\\
\smallskip
\nabla_\xx(\rho \vv\otimes\vv)+\nabla_\xx p=0,\\
\nabla_\xx\cdot(\rho(E+\frac{p}{\rho})\vv)=0,
  \end{array}%
\right.
\end{equation}
where $\vv$ is the velocity,
$E$ is the total energy,
and the constitutive relations among the thermodynamical
variables $\rho, p, e, \theta$, and $S$ are  determined by
(\ref{equilibrium})--(\ref{1.5}).

For the barotropic (isentropic or isothermal) case, $p=p(\rho)$
is determined by \eqref{1.7} with $\gamma\ge 1$,
and then the first $d+1$ equations in \eqref{eq1.1} form
a self-contained system, the Euler system
for steady barotropic fluids.

System (\ref{eq1.1}) governing a supersonic
flow (i.e., $|\vv|^2>c^2$)
has all real eigenvalues and is hyperbolic,
while system (\ref{eq1.1}) governing a subsonic
flow (i.e., $|\vv|^2<c^2$) has complex eigenvalues and
is elliptic-hyperbolic mixed and composite.

\subsection{Wedge Problems Involving Supersonic Shock-Fronts}
The analysis of 2-D steady supersonic
flows past wedges whose vertex angles
are less than the critical angle can date back to the 1940s
since the stability of such flows is fundamental
in applications (cf. \cite{CFr,Wh}).
Local
solutions around the wedge vertex were first constructed in
Gu \cite{Gu}, Li \cite{Lit}, and Schaeffer \cite{Schaeffer}.
Also see Zhang \cite{Zh1} and the references cited therein
for global
potential solutions
when the wedges are a small perturbation of the straight-sided wedge.
For the wedge problem, when the vertex angle is
suitably large, the flow contains a large shock-front
and, for this case, the full Euler equations
(\ref{eq1.1}) are required to describe the physical flow.
When a wedge is straight and its vertex angle is less than
the critical angle $\omega_{c}$,
there exists a supersonic shock-front emanating from the wedge
vertex so that the constant states on both sides
of the shock are supersonic; the critical angle condition
is necessary and sufficient for the existence of the supersonic
shock
(also see \cite{CZZ:1,CFr}).

\smallskip
Consider 2-D steady supersonic
Euler flows past a 2-D Lipschitz curved wedge
$|x_2|\le g(x_1),  x_1>0$, with $g\in Lip(\R_+)$ and $g(0)=0$,
whose vertex angle $\omega_0:=\arctan (g'(0+))$
is less than the critical angle $\omega_{c}$,
along which $TV\{g'(\cdot); \, \R_+\}\le \varepsilon$
for some constant $\varepsilon>0$.
Denote
$$
\Omega:=\{ \xx\,: \, x_2> g(x_1),\, x_1\ge 0 \}, \quad
 \Gamma:= \{ \xx\,:\, x_2= g(x_1),\, x_1\ge 0 \},
$$
and ${\bf n}(x_1\pm) =\frac{(-g'(x_1\pm),1)}{\sqrt{(g'(x_1\pm))^2+1}}$
are the outer normal vectors to $\Gamma$ at points $x_1\pm$, respectively.
The uniform upstream flow $\uu_-=(\rho_-, \rho_-u_-, 0, \rho_-E_-)$
satisfies
$
u_- > c(\rho_-, S_-)
$
so that a strong supersonic shock-front emanates from the wedge
vertex.
Since  the problem is symmetric with respect to the $x_2$-axis,
the wedge problem can be formulated into the
following problem of initial-boundary value type for
system (\ref{eq1.1}) in $\Omega$:
\begin{eqnarray}
&\mbox{\it Cauchy Condition:} \qquad\qquad\quad
   &\u|_{x_1=0} = \u_-; \qquad\qquad\qquad\label{IC} \\
&\mbox{\it Boundary Condition}: \qquad\qquad
  &\vv \cdot {\bf n}=0\qquad\mbox{ on }\, \Gamma.\qquad\qquad\qquad \label{BC}
\end{eqnarray}

In Chen-Zhang-Zhu \cite{CZZ:1},
it has been established that
there exist $\e_0>0$ and $C>0$ such that,
when $\e\le \e_0$,
there exists a pair of functions
$$
\u=(\rho,\vv, \rho E)\in BV(\R;\R^2\times\R_+\times\R_+),\qquad
\sigma\in BV(\R_+;\R)
$$
with $\chi=\int_{0}^{x_1} \sigma(s)ds \in Lip (\R_+;\R_+)$
such that
$\u$ is a global entropy solution of problem
\eqref{eq1.1}--\eqref{BC} in $\Omega$ with
\begin{eqnarray}
&&TV\{ \u(x_1,\cdot)\,:\, [g(x_1), -\infty)\} \le C\, TV\{g'(\cdot)\}
\quad \mbox{for every}\,\, \, x_1 \in \R_+,
\end{eqnarray}
and
the strong shock-front $x_2=\chi(x_1)$
emanating from the wedge
vertex is nonlinearly stable in structure.
Furthermore, the global $L^1$-stability of entropy solutions with respect
to the incoming flow at $x_1=0$ in $L^1$
has also been established in Chen-Li \cite{ChenLi}.
This asserts that any supersonic shock-front
for the wedge problem is nonlinearly $L^1$--stable with respect to
the $BV$
perturbation
of the incoming flow and the wedge boundary.

In order to achieve this, we
have first developed an adaptation of the Glimm scheme whose mesh grids
are designed to follow the slope of the Lipschitz wedge boundary
so that
the lateral Riemann building blocks contain only
one shock or rarefaction wave emanating from the mesh points on
the boundary.
Such a design makes the $BV$ estimates more convenient for
the Glimm approximate solutions.
Then careful interaction estimates have been made.
One of the essential estimates is the estimate of the strength $\delta_1$
of the reflected 1-waves in the interaction between
the 4-strong shock-front and weak waves
$(\alpha_1,\beta_2,\beta_3,\beta_4)$, that is,
$$
\delta_1=\alpha_1+K_{s1}\beta_4 + O(1)|\alpha_1|(|\beta_2|+|\beta_3|)
\qquad\mbox{with}\,\,\,\, |K_{s1}|<1.
$$
The second essential estimate
is the interaction estimate between
the wedge boundary and weak waves.
Based on the construction of the approximation solutions
and interaction estimates,
we have successfully identified a Glimm-type
functional to incorporate the curved wedge boundary and the strong
shock-front naturally and to trace the interactions not only
between the wedge boundary and weak waves, but also
between the strong shock-front and weak waves.
With the aid of the important fact that $|K_{s1}|<1$, we
have showed that the identified Glimm functional
monotonically decreases in the flow direction.
Another essential estimate is to trace the approximate strong
shock-fronts in order to establish the nonlinear stability and asymptotic
behavior of the strong shock-front emanating from the wedge vertex
under the $BV$ wedge perturbation.

\smallskip
For the 3-D {\bf cone problem}, the nonlinear
stability of a self-similar 3-D full gas flow past an
infinite cone is another important problem.
See Lien-Liu \cite{LL} for the cones with small vertex
angle.
Also see \cite{ChS4,CXY} for the construction of piecewise smooth
potential flows under smooth perturbation of the straight-sided cone.

\smallskip
In Elling-Liu \cite{EllingLiu08}, an evidence has been provided
that the steady supersonic weak shock solution
is dynamically stable, in the sense that it describes the
long-time behavior of an unsteady flow.

\subsection{Stability of Supersonic Vortex Sheets}
Another natural problem is the stability of
supersonic vortex sheets above the Lipschitz wall
$x_2=g(x_1), x_1\ge 0,$ with
$$
g\in Lip(\R_+;\R),\,\,\, g(0)=g'(0+)=0,\,\,\,
\lim\limits_{x_1\to\infty}\arctan(g^{\prime}(x_1+))=0,
$$
and $g'\in BV(\R_+;\R)$ such that
$TV\{g'(\cdot)\}\le \varepsilon$ for some $\varepsilon>0$.
Denote again
$
\Omega=\{\xx\, :\, x_2> g(x_1),\, x_1\ge 0\}$ and
$\Gm= \{\xx\, :\,  x_2= g(x_1),\, x_1\ge 0 \}$.
The upstream flow consists of one supersonic straight
vortex sheet $x_2=y_0>0$ and two constant vectors
$\u_0=(\rho_0,\rho_0u_0, 0, \rho_0E_0)$ when $x_2>y_0>0$ and
$\u_1=(\rho_1, \rho_1u_1, 0, \rho_1E_1)$ when $0<x_2<y_0$
satisfying
$u_1>u_0>0$ and $u_i >c(\rho_i, S_i)$ for $i=0,1$.
Then the vortex sheet problem can be formulated
into the following problem of initial-boundary value type for
system (\ref{eq1.1}):
\begin{eqnarray}
&\mbox{\it Cauchy Condition:} \qquad\qquad\quad\,\,\,
   &\u|_{x_1=0} =\left\{\begin{array}{ll} \u_0, \quad 0<x_2<y_0,\\
                 \u_1, \quad  x_2>y_0;
          \end{array}\right. \qquad  \label{IC-2} \\
&\mbox{\it  Slip Boundary Condition}: \qquad  &\vv \cdot {\bf n}=0\qquad\mbox{ on }\, \Gamma.\qquad \label{BC-2}
\end{eqnarray}

It has been proved that
steady supersonic vortex sheets, as time-asymptotics,
are stable in structure globally, even under the $BV$ perturbation of
the Lipschitz walls in Chen-Zhang-Zhu \cite{CZZ:2}.
The result
indicates that the strong supersonic
vortex sheets are nonlinearly stable in structure globally
under the $BV$ perturbation of the Lipschitz wall,
although there may be weak shocks and supersonic vortex sheets away
from the strong vortex sheet.
In order to establish this theorem,
as in \S 4.1, we first developed an adaption of the Glimm scheme
whose mesh grids
are designed to follow the slope of the Lipschitz boundary.
For this case, one of the essential estimates is
the estimate of the strength $\delta_1$ of the reflected 1-wave
in the interaction between the 4-weak wave $\alpha_4$
and the strong vortex sheet from below is less
than one: $\delta_1=K_{01}\alpha_4$ and $|K_{01}|<1$.
The second new essential estimate is the estimate of the strength $\delta_4$
of the reflected 4-wave in the interaction between
the 1-weak wave $\beta_1$  and the strong vortex sheet from above
is also less than one: $\delta_4=K_{11}\beta_1$ and $|K_{11}|<1$.
Another essential estimate is to trace the approximate
supersonic vortex sheets
under the $BV$ boundary
perturbation. For more details, see Chen-Zhang-Zhu \cite{CZZ:2}.
The nonlinear $L^1$-stability of entropy solutions with respect to
the incoming flow at $x_1=0$ is under current investigation.

\section{Multidimensional Steady Transonic Problems}

In this section we discuss another important class of
M-D steady problems: transonic problems.
In the last decade, a program has been initiated on the existence and
stability of M-D transonic
shock-fronts, and some new analytical approaches including
techniques, methods, and ideas have been developed.
For clear presentation, we focus mainly on the celebrated steady
potential flow equation of aerodynamics for the velocity potential
$\varphi: \Omega\subset\bR^d\rightarrow\bR$, which is a second-order
nonlinear PDE of mixed elliptic-hyperbolic type:
\begin{equation}\label{PotenEulerCompres}
\nabla_\xx\cdot (\rho(|\nabla_\xx\varphi|^2)\nabla_\xx\varphi)=0,
\qquad \xx\in\Omega\subset\bR^d,
\end{equation}
where the density $\rho(q^2)$ by scaling is
\begin{equation}\label{PotenEulerCompres-1}
\rho(q^2)=\big(1-\frac{\gamma-1}{2}q^2\big)^{\frac{1}{\gamma-1}}
\end{equation}
with adiabatic exponent $\gamma>1$.
Equation (\ref{PotenEulerCompres})
is elliptic at $\nabla_\xx\varphi$ with $|\nabla_\xx\varphi|=q$ if
$
\rho(q^2)+ 2q^2\rho'(q^2)>0,
$
which is equivalent to
$$
q< q_*:=\sqrt{2/(\gamma+1)},
$$
i.e., the flow is subsonic.
Equation (\ref{PotenEulerCompres}) is hyperbolic if
$
\rho(q^2)+ 2q^2\rho'(q^2)<0,
$
i.e., $q>q_*$, that is, the flow is supersonic.

Let $\Omega^+$ and $\Omega^-$ be open subsets of $\Omega$
such that
$\Omega^+\cap \Omega^-=\emptyset$,
$\overline{\Omega^+}\cup \overline{\Omega^-}=\overline\Omega$,
and $\bS=\partial\Omega^+\cap\Omega$.
Let $\varphi\in C^{0,1}(\Omega)\cap C^1(\overline{\Omega^\pm})$ be a weak solution of
(\ref{PotenEulerCompres}), which
satisfies $\displaystyle |\nabla_\xx\varphi|\leq q_{cav}:=\sqrt{2/(\gamma-1)}$,
so that $\nabla_\xx\varphi$ experiences a jump across $\bS$ that is
a $(d-1)$-D smooth surface.
Set $\varphi^\pm=\varphi|_{\Omega^\pm}$.
Then $\varphi$
satisfies  the
following Rankine-Hugoniot conditions
on $\bS$:
\begin{eqnarray}
&\varphi^+=\varphi^-,
\label{FBCondition-0}\\
&\rho(|\nabla_\xx\varphi^+|^2) \nabla_\xx\varphi^+\cdot {\bf n}
 =\rho(|\nabla_\xx\varphi^-|^2)
  \nabla_\xx\varphi^-\cdot {\bf n}.
  \label{FBCondition-1}
\end{eqnarray}

Suppose that $\varphi\in C^1(\overline{\Omega^\pm})$
is a weak solution satisfying
\begin{equation}\label{TnasonicCond}
|\nabla_\xx\varphi|<q_*
\,\, \mbox{in }\, \Omega^+,\;\;
\quad |\nabla_\xx\varphi|>q_* \,\, \mbox{in }\;\Omega^-,
\;\;\quad \nabla_\xx\varphi^\pm\cdot {\bf n}|_\bS>0.
\end{equation}
Then $\varphi$ is a {\it transonic shock-front solution}
with {\it transonic shock-front} $\bS$ dividing
$\Omega$ into the {\em subsonic region} $\Omega^+$
and the {\em supersonic region} $\Omega^-$ and
satisfying the physical entropy
condition (see Courant-Friedrichs \cite{CFr}):
\begin{equation}\label{entropy}
\rho(|\nabla_\xx\varphi^-|^2)<\rho(|\nabla_\xx\varphi^+|^2)
\qquad \mbox{along}\,\, \bS.
\end{equation}
Note that equation (\ref{PotenEulerCompres})
is elliptic in the subsonic region and
hyperbolic in the supersonic region.

\smallskip
As an example, let $(x_1,\xpr)$ be the coordinates
in $\bR^d$, where $x_1\in \bR$
and $\xpr=(x_2,\dots,x_{d})\in\bR^{d-1}$.
Fix $\VV_0\in \bR^d$, and let
$\varphi_0(\xx)\defd \VV_0\cdot\xx$, $\xx\in \bR^d$.
If $\displaystyle |\VV_0|\in (0, q_*)$
(resp. $\displaystyle |\VV_0|\in (q_*, q_{cav}) $),
then $\varphi_0(\xx)$ is a subsonic (resp. supersonic)
solution in $\bR^d$, and $\VV_0=\nabla_\xx\varphi_0$ is its velocity.

\smallskip
Let $q^-_0> 0$ and $\Vtan_0\in \bR^{d-1}$ be such that
the vector $\Vpl_0 \defd (q^-_0,\Vtan_0)$ satisfies
$|\Vpl_0|> q_*$. Then
there exists a unique $q^+_0>0$ such that
\begin{equation} \label{pPlusMinCondit}
\big(1-\frac{\gamma-1}{2}
(|q^+_0|^2+|\Vtan_0|^2)\big)^{\frac{1}{\gamma-1}}q^+_0
=\big(1- \frac{\gamma-1}{2}
(|q^-_0|^2+|\Vtan_0|^2)\big)^{\frac{1}{\gamma-1}}q^-_0.
\end{equation}
The entropy condition (\ref{entropy}) implies $q^+_0< q^-_0$.
By denoting $\Vmn_0\defd (q^+_0, \Vtan_0)$ and defining the functions
$\varphi_0^\pm(\xx)\defd \Vpm_0\cdot \xx$ on $\bR^d$,
then $\varphi_0^+$ (resp. $\varphi_0^-$) is a subsonic (resp. supersonic)
solution. Furthermore,
from (\ref{FBCondition-1}) and (\ref{pPlusMinCondit}),
the function
\begin{equation} \label{linearTransonicSol}
\varphi_0(\xx)\defd \min (\varphi_0^-(\xx), \varphi_0^+(\xx))
=
V^\pm_0\cdot\xx, \,\,\,\,\,\, \xx\in\Omega^\pm_0\defd \{\pm x_1>0\},
\end{equation}
is a plane transonic shock-front solution in $\bR^d$,
$\Omega^-_0$ and $\Omega^+_0$ are respectively its supersonic and
subsonic regions, and $\bS=\{x_1=0\}$ is a transonic shock-front.
Note that, if $\Vtan_0=0$, the velocities $\Vpm_0$ are
orthogonal to the shock-front $\bS$ and, if $\Vtan_0\ne 0$,
the velocities are not orthogonal to $\bS$.

\subsection{Transonic Shock-Front Problems in $\R^d$}
Consider M-D perturbations of
the uniform transonic shock-front solution
(\ref{linearTransonicSol}) in $\bR^d$
with $d\ge 3$.
Since it suffices to specify the supersonic
perturbation $\varphi^-$ only in a neighborhood of the unperturbed
shock-front $\{x_1=0\}$,
we introduce domains
$
\Omega:=(-1,\infty)\times\bR^{d-1},
\Omega_1:=(-1, 1)\times\bR^{d-1}.
$
Note that we expect the subsonic region
 $\Omega^+$ to be close to the half-space
$\Omega^+_0=\{x_1>0\}$.

\medskip
\noindent
{\bf Problem 5.1}. Given a supersonic solution $\varphi^-(\xx)$ of (\ref{PotenEulerCompres})
in $\Omega_1$,
find a transonic shock-front solution $\varphi(\xx)$ in $\Omega$
such that
$$
\Omega^-\subset\Omega_1, \qquad
\varphi(\xx)=\varphi^-(\xx) \,\,\,\,\,\mbox{ in}\,\, \Omega^-,
$$
where $\Omega^-\defd \Omega\setminus\Omega^+$ and
$\Omega^+\defd \{\xx\in\Omega \;: \;|\nabla_\xx\varphi(\xx)|<q_*\}$,
and
\begin{eqnarray}
\label{bdryConditions}
&&\varphi=\varphi^-, \quad
\partial_{x_1}\varphi=\partial_{x_1}\varphi^-
\qquad \mbox{on }  \{x_1=-1\}, \\
\label{stateAtInfty}
&&\lim_{R\to\infty}\|\varphi-\varphi^+_0\|_{C^1(\Omega^+\setminus B_R(0))}=0.
\end{eqnarray}

Condition {\rm (\ref{bdryConditions})} determines that the solution
has supersonic upstream,
while condition {\rm (\ref{stateAtInfty})} especially determines that the
uniform velocity state at infinity in the
downstream direction is equal to the unperturbed downstream velocity
state. The additional requirement
in {\rm (\ref{stateAtInfty})} that $\varphi\to \varphi^+_0$ at infinity
within $\Omega^+$ fixes the position of shock-front at infinity.
This allows us to determine the solution of {\it Problem 5.1} uniquely.
In Chen-Feldman \cite{CFe1}, we have employed the free boundary approach
first developed for the potential flow equation
(\ref{PotenEulerCompres}) to solve the stability problem, {\it Problem 5.1}.
This existence result can be extended  to the case that the regularity
of the steady perturbation $\varphi^-$ is only $C^{1,1}$.
It would be interesting to establish similar results for the full
Euler equations \eqref{Euler1}.

\subsection{Nozzle Problems Involving Transonic Shock-Fronts}
We now consider M-D transonic shock-fronts in
the following infinite nozzle $\Noz$
with arbitrary smooth cross-sections:
$\Noz=\CylToNoz(\Csect \times\R)\cap\{x_1> -1\}$,
where
$\Csect\subset\bR^{d-1}$ is an open bounded connected set with
a smooth boundary,
and $\CylToNoz: \bR^d\to \bR^d$
is a smooth map, which is close to the identity map. For simplicity,
we assume that
$\partial\Csect$ is in $\displaystyle C^{[\frac{d}{2}]+3,\alpha}$
and
$\|\CylToNoz-I\|_{[\frac{d}{2}]+3, \alpha, \bR^d} \leq \epsP$
for some $\alpha\in(0, 1)$ and small $\sigma>0$,
where $[s]$ is the integer part of $s$,
$I: \bR^d\to \bR^d$ is the identity map, and
$\partial_l\Noz\defd \CylToNoz(\R\times\partial\Csect)
\cap\{x_1>-1\}$.
For concreteness, we also assume that
there exists $L>1$
such that
$\CylToNoz(\xx)=\xx$ for any
$\xx=(x_1,\xx')$ with $x_1>L$.

In the 2-D case,
$
\Noz=\{(x_1, x_2)\; : \; x_1>-1, \; b^-(x_2)< x_2< b^+(x_2)\},
$
where
$\|b^\pm-b^\pm_\infty\|_{4,\alpha,\bR}\leq \epsP$ and
$b^\pm\equiv b^\pm_\infty$ on $[L, \infty)$ for some
constants $b^\pm_\infty$ satisfying $b^+_\infty>b^-_\infty$.
For the M-D case, the geometry of
the nozzles is much richer.
Note that our setup implies that
$
\partial\Noz=\overline{\partial_o\Noz}\cup\partial_l\Noz
$
with
\begin{eqnarray*}
&&\partial_l\Noz\defd \CylToNoz[(-\infty, \infty)\times\partial\Csect]
\cap\{(x_1,\xx')\; : \;x_1>-1\},
\\
&&\partial_o\Noz\defd \CylToNoz((-\infty,\infty)\times\Csect)
 \cap\{(x_1,\xx')\; : \;x_1=-1\}.
\end{eqnarray*}
Then our transonic nozzle problem can be
formulated as

\medskip
\noindent
{\bf Problem 5.2: Transonic Nozzle Problem}.
Given the supersonic upstream flow at the entrance $\partial_o\Omega$:
\begin{equation}
\varphi=\varphi^-_e, \,\,\,\,\varphi_{x_1}=\psi^-_e
  \qquad \mbox{on }  \partial_o\Noz,
\label{bdryConditions_inf_1}
\end{equation}
the slip boundary condition on the nozzle boundary $\partial_l\Noz$:
\begin{equation}
\nabla_\xx\varphi\cdot {\bf n}=0 \qquad \mbox{on } \mbox{$\partial_l\Noz$},
\label{bdryConditions_inf_2}
\end{equation}
and the uniform subsonic flow condition at the infinite
exit $x_1=\infty$:
\begin{equation}
\|\varphi(\cdot)- q_\infty x_1\|_{C^1(\Noz\cap\{x_1>R\})}\to 0
\qquad\, \mbox{as } R\to\infty \,\,
\mbox{for some } q_\infty\in (0, q_*),
\label{bdryConditionsatinfinity}
\end{equation}
find a M-D transonic flow $\varphi$
of problem (\ref{PotenEulerCompres}) and
(\ref{bdryConditions_inf_1})--(\ref{bdryConditionsatinfinity})
in $\Omega$.

\medskip
The standard local existence theory of smooth solutions for the
initial-boundary value problem
(\ref{bdryConditions_inf_1})--(\ref{bdryConditions_inf_2})
for second-order quasilinear hyperbolic equations
implies that, as $\sigma$ is sufficiently small,
there exists a supersonic solution $\varphi^-$ of
(\ref{PotenEulerCompres}) in
$
\Omega_2:=\{-1\le x_1\le 1\},
$
which is a $C^{l+1}$ perturbation of $\varphi_0^-=q_0^-x_1$:
For any $\alpha\in (0,1]$,
\begin{equation}\label{smallPert_inf-I}
\|\varphi^- - \varphi_0^-\|_{l, \alpha,\Noz_2}\le C_0\epsP,
\qquad l=1,2,
\end{equation}
for some constant $C_0>0$, and satisfies
$\nabla\varphi^-\cdot {\bf n}=0$
on $\partial_l \Noz_2$,
provided that
$(\varphi^-_e,\psi^-_e)$ on
$\partial_o\Omega$ satisfying
\begin{equation}\label{smallPert_inf-Cauchy}
\|\varphi_e^--q_0^-x_1\|_{H^{s+l}}+\|\psi_e^--q_0^-\|_{H^{s+l-1}}
\le \sigma, \qquad l=1,2,
\end{equation}
for some integer $s>{d}/{2}+1$ and the compatibility conditions up to
order $s+1$, where the norm $\|\cdot\|_{H^s}$ is the Sobolev
norm with $H^s=W^{s,2}$.

\smallskip
In Chen-Feldman \cite{CFe3}, {\it Problem 5.2} has been solved.
More precisely, let $\displaystyle q^-_0\in (q_*, q_{cav})$
and $\displaystyle q^+_0\in (0, q_*)$
satisfy (\ref{pPlusMinCondit}),
and let $\varphi_0$ be the transonic shock-front solution
(\ref{linearTransonicSol}) with $\Vtan=0$.
Then there exist $\epsP_0>0$ and $C$,
depending only on $d$, $\alpha$, $\gamma$, $q^-_0$, $\Csect$, and $L$,
such that,
for every $\epsP \in (0,\epsP_0)$, any  map $\CylToNoz$ introduced above,
and any supersonic upstream flow $(\varphi^-_e,\psi^-_e)$ on
$\partial_o\Omega$ satisfying (\ref{smallPert_inf-Cauchy}) with $l=1$,
there exists a solution
$\varphi\in C^{0,1}(\Noz) \cap C^{1,\alpha}(\overline{\Omega^+})$
satisfying
\begin{eqnarray*}
&&\Omega^+(\varphi)=\{x_1> f(\xpr)\}, \qquad
\Omega^-(\varphi)=\{x_1< f(\xpr)\}, \label{subsup-regions}\\
&&\|\varphi - q^\pm x_1\|_{1,\alpha,\Omega^\pm}
\le C\sigma, \qquad
\|\varphi- q_\infty x_1\|_{C^1(\Noz\cap\{x_1>R\})} \to 0
\,\, \mbox{as $R\to\infty$},
\end{eqnarray*}
where $q_\infty\in (0, q_*)$ is the unique solution of the equation
\begin{equation}\label{DefFBPLinearGrowth}
\rho((q^+)^2)q^+=Q^+:=
\frac{1}{|\Csect|}
\int_{\partial_o\Noz}\rho(|\nabla_{\xx'}\varphi_e^-|^2+(\psi_e^-)^2)
\psi_e^-\,d{{\mathcal H}^{d-1}}.
\end{equation}
In addition, if the supersonic uniform
flow $(\varphi_e^-,\psi_e^-)$ on
$\partial_o\Omega$ satisfies
(\ref{smallPert_inf-Cauchy}) with $l=2$, then
$\varphi\in C^{2,\alpha}(\overline{\Omega^+})$ with
$
\|\varphi- q^+x_1\|_{2,\alpha,\Omega_+}\le C\sigma,
$
and the solution with a transonic shock-front
is unique and stable with respect to the nozzle boundary
and the smooth supersonic upstream flow at the entrance.
The techniques have been extended to solving
the nozzle problem for the 2-D full Euler
equations first in Chen-Chen-Feldman \cite{CCF1}.

\smallskip
In the previous setting, the location of the transonic shock-front
in the solution is not unique in general if we prescribe only
the pressure at the nozzle exit, since the flat shock-front
between uniform states can be translated
along the flat nozzle which does not change the flow parameters
at the entrance and the pressure at the exit.
The analysis of the relation among the
unique shock-front location, the flow parameters, and the geometry of
the diverging nozzle has been made for various Euler systems in
different dimensions;
see Chen \cite{SxChen} and Lin-Yuan \cite{LYu} for the 2-D
full Euler equations,
Li-Xin-Yin \cite{LXY} for the 2-D and 3-D
axisymmetric Euler systems with narrow divergent nozzles,
and Bae-Feldman \cite{BaeF} for perturbed diverging cone-shaped nozzle
of arbitrary cross-section for the non-isentropic potential flow
system in any dimension.

\smallskip
A major open problem is the {\it physical de Laval nozzle problem}:
Consider a nozzle $\Omega$ which is flat (i.e., $\Psi=I$) between
$-2\le x_1\le -1$ and $1\le x_1\le \infty$, and has some special
geometry between $-1\le x_1\le 1$ to make the part $\{-1\le x_1\le 0\}$
convergent and the part $\{0\le x_1\le 1\}$ divergent.
Given certain incoming subsonic flow at $x_1=-2$, find a transonic flow
containing transonic shock-fronts in the nozzle such that the downstream
at $x_1=\infty$ is subsonic.
A related reference is Kuz'min's book \cite{Kuzmin}.

\subsection{Wedge/Cone Problems Involving Transonic Shock-Fronts}
The existence and stability of
transonic flows past infinite wedges or cones
are further long-standing open transonic problems.
Some progress has been made for the wedge case in 2-D
steady flow in Chen-Fang \cite{ChenFang1}, Fang \cite{Fang},
and Chen-Chen-Feldman \cite{CCF2}.
In \cite{Fang}, it was proved that the
transonic shock-front is conditionally stable under perturbation of the
upstream flow and/or the wedge boundary in some weak Sobolev
norms.
In \cite{CCF2}, the
existence and stability of transonic flows past the  curved wedge have been
established in the strong H\"{o}lder norms for the full Euler equations.

Conical flow (i.e. cylindrically symmetric flow
with respect to an axis, say, the $x_1$-axis) in $\R^3$ occurs in many physical
situations (cf. \cite{CFr}).
Unlike the 2-D case, the governing equations for the 3-D conical case have
a singularity at the cone vertex and
the flow past the straight-sided cone is self-similar, but is no
longer piecewise constant.
These have resulted in additional difficulties for the stability problem.
In Chen-Fang \cite{CFang},
we have developed techniques to handle the singular terms in
the equations and the singularity of the solutions.
Our main results indicate that the self-similar transonic
shock-front is conditionally stable with respect to the conical
perturbation of the cone boundary and the upstream flow in
appropriate function spaces.
That is, the transonic shock-front and downstream
flow in our solutions are close to the unperturbed self-similar
transonic shock-front and downstream flow under the conical
perturbation, and
the slope
of the shock-front asymptotically tends to the slope of the
unperturbed self-similar shock-front at infinity.

In order to achieve these results, we have first formulated the stability
problem as a free boundary problem and have then introduced a coordinate
transformation to reduce the free boundary problem into a fixed
boundary value problem for a singular nonlinear elliptic system. We
have developed an iteration scheme that consists of two iteration mappings:
one is for an iteration of approximate transonic shock-fronts; and
the other is for an iteration of the corresponding boundary value
problems for the singular nonlinear systems for given approximate
shock-fronts. To ensure the well-definedness and contraction
property of the iteration mappings, it is essential to establish the
well-posedness for a corresponding singular linearized elliptic
equation, especially the stability with respect to the coefficients
of the equation, and to obtain the estimates of its solutions
reflecting their singularity at the cone vertex and decay at
infinity. The approach is to employ key features of the equation,
introduce appropriate solution spaces, and apply a Fredholm-type
theorem in Maz'ya-Plamenevski\v{\i} \cite{MP} to establish the
existence of solutions by showing the uniqueness in the solution
spaces. Also see Cui-Yin \cite{CuiYin} for related results.

\subsection{Airfoil/Obstacle Problems: Subsonic Flows past
an Airfoil or an Obstacle}

The Euler equations for potential flows
\eqref{PotenEulerCompres}--\eqref{PotenEulerCompres-1}
in $\R^2$ can be rewritten as
\begin{equation}\label{5.1}
\begin{cases}
\del_{x_1}v -\del_{x_2} u = 0,\\
\del_{x_1}(\rho u)+ \del_{x_2}(\rho v) = 0,
\end{cases}
\end{equation}
with $(u,v)=\nabla_\xx\varphi$ and $q=|\nabla_\xx\varphi|$,
where $\rho$ is again determined by the Bernoulli relation
\eqref{PotenEulerCompres}.
An important problem is subsonic flows past an airfoil
or an obstacle.
Shiffman \cite{Sh}, Bers \cite{Bers}, Finn-Gilbarg \cite{FinnGilbarg},
and Dong \cite{Dong} studied
subsonic (elliptic) solutions of (\ref{PotenEulerCompres}) outside
an obstacle when the upstream flows are sufficiently subsonic;
also see Chen-Dafermos-Slemrod-Wang \cite{CDSW} via compensated compactness
argument.
Morawetz in \cite{M3} first showed that the
flows of (\ref{5.1}) past an obstacle  may
contain transonic shocks in general.
A further problem is to construct global entropy solutions
of the airfoil problem (see \cite{Mora1,GM1}).

In Chen-Slemrod-Wang \cite{CSW1},
we have
introduced the usual flow angle
$\theta = \tan^{-1}(\frac{v}{u})$ and written the irrotationality
and mass conservation equation as an artificially viscous problem:
\begin{equation} \label{5.15}
\begin{cases}
   \partial_{x_1} v^\e -\partial_{x_2} u^\e=\varepsilon\Delta\theta^\e ,  \\
   \partial_{x_1}(\rho^\e u^\e) +\partial_{x_2}(\rho^\e v^\e)
   =\varepsilon\nabla\cdot(\sigma(\rho^\e)\nabla \rho^\e),
\end{cases}
\end{equation}
where $\sigma (\rho)$ is suitably chosen, and appropriate boundary
conditions are imposed for this regularized ``viscous'' problem. The
crucial new discovery is that a uniformly $L^\infty$ bound in
$q^\varepsilon$ can be obtained when $1\leq\gamma <3$ which uniformly
prevents cavitation.  However, in this formulation,
a uniform bound in the flow angle
$\theta^\varepsilon$ and a uniform lower bound in $q^\varepsilon$
in any fixed region disjoint
from the profile must be assumed apriori.
By making further careful energy estimates,
Morawetz's argument \cite{Mora1} then applies, and the strong
convergence in $L^1_{\rm loc} (\Omega)$ of our approximating
sequence is achieved.

An important open problem is whether two remaining conditions
on $(q^\varepsilon, \theta^\varepsilon)$ can be removed.
See Chen-Slemrod-Wang \cite{CSW1} for more details.

\subsection{Nonlinear Approaches}
We now discuss several nonlinear approaches to deal with
steady transonic problems.

\subsubsection{Free Boundary Approaches}
We first describe two of the free boundary approaches
for {\it Problems 5.1--5.2}, developed in
\cite{CFe1,CFe3}.

{\bf Free Boundary Problems}.
The transonic shock-front problems can be formulated
into a one-phase free boundary problem for a nonlinear
elliptic PDE:
Given $\varphi^-\in C^{1,\alpha}(\overline\Omega)$,
find a function $\varphi$ that is continuous in $\Omega$
and satisfies
\begin{equation} \label{transonicInequalitiesPert}
\varphi \leq \varphi^- \qquad \hbox{in}\,\, \bar{\Omega},
\end{equation}
equation (\ref{PotenEulerCompres}), the ellipticity
condition in the non-coincidence set
$\Omega^+=\{\varphi<\varphi^-\}$,
the free boundary condition (\ref{FBCondition-1})
on the boundary $\mathcal{S}=\partial\Omega^+\cap \Omega$,
as well as the prescribed conditions on the fixed
boundary $\partial\Omega$ and at infinity.
These conditions are different in different problems,
for example, conditions
(\ref{bdryConditions})--(\ref{stateAtInfty})
for {\it Problem 5.1} and
(\ref{bdryConditions_inf_1})--(\ref{bdryConditionsatinfinity})
for {\it Problem 5.2}.

The free boundary is the location of the shock-front, and the free
boundary conditions (\ref{FBCondition-0})--(\ref{FBCondition-1})
are the Rankine-Hugoniot conditions.
Note that condition (\ref{transonicInequalitiesPert}) is motivated
by the similar property (\ref{linearTransonicSol})
of unperturbed shock-fronts;
and (\ref{transonicInequalitiesPert}), locally on the shock-front,
is equivalent to the entropy condition (\ref{entropy}).
Condition  (\ref{transonicInequalitiesPert})
transforms the transonic shock-front problem, in which the
subsonic region $\Omega^+$ is determined by the gradient condition
$|\nabla_\xx\varphi(\xx)|<q_*$, into a free boundary problem
in which $\Omega^+$ is the non-coincidence set.
In order to solve this free boundary problem,
equation (\ref{PotenEulerCompres}) is modified to be
uniformly elliptic and then
the free boundary condition (\ref{FBCondition-1}) is
correspondingly modified.
The problem becomes a one-phase free boundary problem
for the uniformly elliptic equation which we can solve.
Since
$\varphi^-$ is a small $C^{1, \alpha}$--perturbation of
$\varphi^-_0$,
the solution $\varphi$ of the free boundary problem is shown
to be a small $C^{1, \alpha}$--perturbation of the given subsonic
shock-front solution $\varphi^+_0$ in $\Omega^+$.
In particular, the gradient estimate implies that $\varphi$
in fact satisfies the original free boundary problem, hence
the transonic shock-front problem, {\it Problem 5.1}
({\it Problem 5.2}, respectively).

The modified free boundary problem
does not directly fit into the variational framework
of Alt-Caffarelli \cite{AC} and Alt-Caffarelli-Friedman \cite{ACF},
and the regularization framework of
Berestycki-Caffarelli-Nirenberg \cite{BCN1}.
Also, the nonlinearity of the free boundary problem makes it difficult
to apply the Harnack inequality approach of Caffarelli \cite{Ca}.
In particular, a boundary comparison principle for positive solutions
of nonlinear elliptic equations in Lipschitz domains is not available yet
for the equations that are not homogeneous with respect to
$(\nabla^2_\xx u, \nabla_\xx u, u)$, which is our case.

\smallskip
{\bf Iteration Approach}.
The first approach we developed in
Chen-Feldman \cite{CFe1,CFe3} for the steady potential flow equation
is an iteration scheme based on the
non-degeneracy of the free boundary condition:
the jump of the normal derivative of a solution across the free boundary
has a strictly positive lower bound.
Our iteration process is as follows. Suppose that the domain $\Omega^+_k$
is given so that
$S_k\defd \partial \Omega^+_k \setminus\partial \Omega$ is
$C^{1, \alpha}$. Consider the oblique derivative problem
in $\Omega^+_k$ obtained by rewriting the (modified)
equation (\ref{PotenEulerCompres}) and free boundary condition
(\ref{FBCondition-1}) in terms of the function
$u\defd \varphi-\varphi^+_0$. Then the problem has the
following form:
\begin{equation}\label{generalFBP}
\begin{array} {ll}
\displaystyle
\divg_\xx \AAA(\xx, \nabla_\xx u)=F(\xx)\;
 \;\;& \mbox{ in }\; \Omega^+_k\defd \{u>0\},\\
\displaystyle \AAA(\xx, \nabla_\xx u)\cdot {\bf n}
   = G(\xx, {\bf n})\; & \mbox{ on }\;
  \bS \defd \partial\Omega^+_k\setminus \partial \Omega,
\end{array}
\end{equation}
plus the fixed boundary conditions
on $\partial\Omega^+_k\cap\partial\Omega$ and the conditions at infinity.
The equation is quasilinear, uniformly elliptic,
$\AAA(\xx, 0)\equiv 0$,
while
$G(\xx,{\bf n})$ has a certain
structure. Let $u_k\in C^{1, \alpha}(\overline{\Omega_k^+})$
be the solution of (\ref{generalFBP}).
Then $\|u_k\|_{1,\alpha,\Omega^+_k}$ is estimated to be small
if the perturbation is small, where appropriate weighted
H\"{o}lder norms are actually needed
in the unbounded domains.
The function $\varphi_k\defd \varphi^+_0+ u_k$
from $\Omega^+_k$ is extended to $\Omega$
so that the $C^{1, \alpha}$--norm
of $\varphi_k-\varphi^+_0$ in $\Omega$ is controlled
by $\|u_k\|_{1,\alpha,\Omega_k^+}$.
For the next step, define
$
\Omega^+_{k+1}\defd\{\xx \in \Omega\,  :\, \varphi_k(\xx)< \varphi^-(\xx)\}.
$
Note that, since $\|\varphi_k-\varphi^+_0\|_{1,\alpha,\Omega}$
and $\|\varphi^--\varphi^-_0\|_{1,\alpha,\Omega}$
are small, we have
$|\nabla_\xx\varphi^-|-|\nabla_\xx\varphi_k|\geq \delta>0$
in $\Omega$,
and this nondegeneracy implies that
$S_{k+1}\defd \partial \Omega^+_{k+1} \setminus\partial \Omega$ is
$C^{1, \alpha}$ and its norm is estimated in terms of the data of
the problem.

The fixed point $\Omega^+$ of this process determines a solution of
the free boundary
problem since the corresponding solution $\varphi$
satisfies $\Omega^+=\{\varphi<\varphi^-\}$
and  the Rankine-Hugoniot condition \eqref{FBCondition-1}
holds on
$\bS\defd \partial\Omega^+ \cap\Omega$.
On the other hand, the elliptic estimates alone are not sufficient
to get the existence of a fixed point, because
the right-hand side of the boundary condition
in problem (\ref{generalFBP})
depends on the unit normal ${\bf n}$ of the free boundary.
One way is to require the orthogonality of the flat
shock-fronts so that
$\rho(|\nabla_\xx\varphi_0^+|^2)
\nabla_\xx\varphi_0^+ = \rho(|\nabla_\xx\varphi_0^-|^2)
\nabla_\xx\varphi_0^-$
in $\Omega$
to obtain better estimates for the iteration and
prove the existence of a fixed point.
For more details,
see Chen-Feldman \cite{CFe1,CFe3}.

\smallskip
{\bf Partial Hodograph Approach}.
The second approach we have developed in \cite{CFe1}
is a partial hodograph procedure,
with which we can handle
the existence and stability of M-D
transonic shock-fronts that are
not nearly orthogonal to the flow direction.
One of the main ingredients in this new approach is to employ a partial
hodograph transform to reduce the free boundary problem
to a conormal boundary value problem
for the corresponding nonlinear
second-order elliptic equation of divergence form
and then to develop
techniques to solve the conormal boundary value problem.
To achieve this,
the strategy is to construct first solutions in
the intersection domains between the physical unbounded domain under
consideration and a series of half
balls with radius $R$, then make uniform estimates in $R$, and
finally send $R\to\infty$.
It requires delicate apriori estimates to achieve this.
A uniform bound in a weighted
$L^\infty$-norm can be achieved
by both employing a comparison principle and identifying
a global function with the same decay rate as the fundamental solution
of the elliptic equation with constant coefficients which controls
the solutions.
Then, by scaling arguments, the uniform estimates can be obtained
in a weighted H\"{o}lder norm for the solutions,
which lead to the existence of a solution in the unbounded domain
with some decay rate at infinity.
For such decaying solutions,
a comparison principle holds, which implies the
uniqueness for the conormal problem.
Finally, by the gradient estimate, the limit function
can be shown to be a solution of
the M-D transonic shock problem,
and then the existence result can be extended to the case that the
regularity of the steady perturbation is only $C^{1,1}$.
We can further prove that the M-D transonic shock-front solution
is stable with respect to the $C^{2,\alpha}$
supersonic perturbation.

\smallskip
When the regularity of the steady perturbation is $C^{3,\alpha}$ or higher:
$
\|\varphi^- - \varphi_0^-\|^{(d-1)}_{3, \alpha,\Omega_1}\le \sigma,
$
we have introduced another simpler approach, the implicit function approach, to deal with
the existence and stability problem in Chen-Feldman \cite{CFe1}.

\smallskip
Other approaches can be found in Canic-Keyfitz-Lieberman
\cite{CanicKeyfitz},
Canic-Keyfitz-Kim \cite{CanicKeyfitzKim},
Chen-Chen-Feldman \cite{CCF1},
Chen \cite{ChS5,SxChen},
Zheng \cite{Zhe2}, and the references cited therein.

\subsubsection{Weak Convergence Approaches}
Recent investigations on conservation laws based on weak convergence
methods
suggested that the
method of compensated compactness be amenable to flows
which exhibit both elliptic and hyperbolic regimes.
In \cite{Mora1} (also see
\cite{Mora3}), Morawetz layed out a program for proving the
existence of the steady transonic flow problem about a bump profile
in the upper half plane (which is equivalent to a symmetric profile
in the whole plane).
It was shown
that, if the key hypotheses of the method of compensated compactness
could be satisfied, now known as a ``compactness framework'' (see
Chen \cite{Chen1}), then indeed there would exist a weak solution
to the problem of flow over a bump which is exhibited by subsonic
and supersonic regimes, i.e., transonic flow.

The compactness framework for system \eqref{5.1} can be formulated
as follows. Let a sequence of functions $\ww^\varepsilon
(\xx)=(u^\varepsilon, v^\varepsilon)(\xx)$ defined on an open set
$\Omega\subset\R^2$ satisfy the following  set of conditions:

\smallskip
\noindent (A.1)\,
 $q^\varepsilon(\xx)=|\ww^\varepsilon (\xx) | \leq q_\sharp
  \;\; \hbox{a.e. in} \;\; \Omega$
for some positive constant $q_\sharp < q_{cav}$;

\smallskip
\noindent (A.2) \, $\nabla_\xx\cdot \QQ_{\pm}(\ww^\varepsilon )$\,\,
are confined in a
compact set in $H^{-1}_{\rm loc}(\Omega)$ for entropy pairs
$\QQ_\pm=(Q_{1\pm}, Q_{2\pm})$,
where $\QQ_{\pm}(\ww^\varepsilon)$ are confined to a bounded set
uniformly in $L_{\rm loc}^\infty
(\Omega)$.

\smallskip
When (A.1) and (A.2) hold,
the Young measure $\nu_{\xx}(\ww), \ww=(u,v)$, determined by the
uniformly bounded sequence of functions $\ww^\varepsilon$
is constrained by the following commutator relation:
\begin{equation}\label{5.14}
 \langle \nu_\xx,\;  Q_{1+} Q_{2-}- Q_{1-}Q_{2+}\rangle
  = \langle \nu_\xx, \; Q_{1+}\rangle \langle \nu_\xx, \; Q_{2-}\rangle
 - \langle \nu_\xx, \; Q_{1-}\rangle
\langle \nu_\xx, \; Q_{2+}\rangle.
\end{equation}
The main point for the compensated compactness framework is to prove
that $\nu_{\xx}$ is a Dirac mass by using entropy pairs, which
implies the compactness of the sequence
$\ww^\varepsilon(\xx)=(u^\varepsilon, v^\varepsilon)(\xx)$ in
$L^1_{\rm loc}(\Omega)$.  In this context,
It is needed Morawetz \cite{Mora1} to presume
the existence of an approximating sequence parameterized by
$\varepsilon$ to their problems satisfying (A.1) and (A.2) so that
they could exploit the commutator identity and obtain the strong
convergence in $L^1_{\rm loc}(\Omega)$ to a weak solution of their
problems.

As it turns out, there is a classical problem where (A.1) and (A.2) hold
trivially, i.e., the sonic limit of subsonic flows.  In that case,
we return to the result by Bers \cite{Bers} and Shiffman
\cite{Sh}, which says that, if the speed at infinity, $q_\infty$,
is less than some $\hat q$, there is a smooth unique solution of
the problem
and ask what
happens as $q_\infty \nearrow \hat q$.  In this case, the flow
develops sonic points and the governing equations become degenerate
elliptic. Thus, if we set $\varepsilon=\hat q-q_\infty$ and examine
a sequence of exact smooth solutions to our system, we see trivially
that (A.1) is satisfied since $|q^\varepsilon|\leq q_*$, and
(A.2) is also satisfied since
$\nabla_\xx\cdot \QQ_\pm (\ww^\varepsilon)=0$
along our solution  sequence.
The effort is in finding entropy pairs which guarantee the
Young measure $\nu_{\xx}$ reduces to a Dirac mass. Ironically, the
original conservation equations of momentum in fact provide two sets
of entropy pairs, while the irrotationality and mass conservation
equations provide another two sets.  This observation has been
explored in detail in Chen-Dafermos-Slemrod-Wang \cite{CDSW}.

What then about the fully transonic problem of flow past an obstacle
or bump where $q_\infty > \hat q$? In Chen-Slemrod-Wang \cite{CSW1},
we have provided some of the ingredients to satisfying (A.1) and
(A.2) as explained in \S 5.4.
On the other hand, (A.2) is easily
obtained from the viscous formulation by using a special entropy
pair of Osher-Hafez-Whitlow \cite{OHW}.  In fact, this entropy pair
is very important:  It guarantees that the inviscid limit of the
above viscous system satisfies a physically meaningful ``entropy''
condition (Theorem 2 of \cite{OHW}). With (A.1)
and (A.2) satisfied, the compensated compactness argument
then applies, to yield the strong
convergence in $L^1_{\rm loc} (\Omega)$ of our approximate solutions.

Some further compensated compactness frameworks have been
developed for solving the weak continuity of
solutions to the Gauss-Codazzi-Ricci equations
in Chen-Slemrod-Wang \cite{CSW2,CSW3}.
Also see Dacorogna \cite{Dac}, Evans \cite{Ev1}, and the references
cited therein for various weak convergence methods and techniques.

\section{Shock Reflection-Diffraction and Self-Similar Solutions}
One of the most challenging problems is to
study solutions with data that give rise to self-similar
solutions (such solutions include
Riemann solutions).
For the Euler equations \eqref{Euler1} for $\xx\in\R^2$, the self-similar
solutions
\begin{equation}\label{self-similar}
(\rho, \vv, p)
=(\rho, \vv, p)(\xi,\eta), \qquad (\xi,\eta)=\xx/{t},
\end{equation}
are determined by
\begin{equation}
\left\{\begin{array}{l}
\del_\xi(\rho U) +\del_\eta(\rho V)=-2\rho,\\
\del_\xi(\rho U^2+p) +\del_\eta(\rho UV)=-3\rho U,\\
\del_\xi(\rho UV) +\del_\eta(\rho V^2 + p)=-3\rho V,\\
\del_\xi(\rho U(E + p/\rho)) +\del_\eta(\rho V(E+ p/\rho)) =-2\rho(E+p/\rho),
\end{array}\right.
\label{Euler1a}
\end{equation}
where $(U,V)=\vv-(\xi, \eta)$ is the pseudo-velocity
and $E=e +\frac{1}{2}(U^2+V^2)$.

The four eigenvalues are
\begin{eqnarray*}
\lambda_0=\frac{V}{U} \,\,\,\mbox{(2-multiplicity)},
\qquad \lambda_\pm=\frac{UV\pm c\sqrt{U^2+V^2-c^2}}{U^2-c^2},
\end{eqnarray*}
where $c=\sqrt{p_\rho(\rho, S)}$ is the sonic speed.
When $U^2+V^2>c^2$, system \eqref{Euler1a} is hyperbolic with four real
eigenvalues and the flow is called pseudo-supersonic, simply called supersonic
without confusion.
When $U^2+V^2<c^2$, system \eqref{Euler1a} is hyperbolic-elliptic
composite type (two repeated eigenvalues are real
and the other two are complex):
two equations are hyperbolic and the other two are elliptic.
The region $U^2+V^2=c^2$ in the $(\xi,\eta)$-plane is simply called
the sonic region in the flow.
In general, system \eqref{Euler1a} is hyperbolic-elliptic mixed
and composite type, and the flow is transonic.
For a bounded solution $(\rho,\vv, p)$,
the flow must be supersonic when $\xi^2+\eta^2\to \infty$.

An important prototype problem
for both practical applications and the theory of
M-D complex wave patterns is the problem
of diffraction of a shock-front which is
incident along an inclined ramp.
When a plane shock-front hits a wedge head on, a self-similar reflected
shock-front moves outward as the original shock-front moves forward
(cf. \cite{BenDor,CC,CFr,GM,Mora2,Serreref,VanDyke}).
Then the problem of shock reflection-diffraction by a wedge can be
formulated as follows:

\medskip
{\bf Problem 6.1 (Initial-boundary value problem)}. {\it Seek a
solution of system \eqref{Euler1} satisfying the initial condition at
$t=0$:
\begin{equation} \label{ibv-1}
(\rho,\vv, p) =\begin{cases}
(\rho_0, 0,0,p_0), &\, \quad |x_2|>x_1\tan\theta_w, x_1>0,\\
(\rho_0, u_1,0,p_1), &\, \quad x_1<0;
\end{cases}
\end{equation}
and the slip boundary condition along the wedge boundary:
\begin{equation}\label{boundary-condition}
\vv\cdot \bnu=0,
\end{equation}
where $\bnu$ is the exterior unit normal to the wedge boundary,
and state $(0)$ and $(1)$ satisfy
\begin{equation}\label{2.10}
u_1=\sqrt{\frac{(p_1-p_0)(\rho_1-\rho_0)}{\rho_0\rho_1}}, \,\,
\frac{p_1}{p_0}=\frac{(\gamma+1)\rho_1-(\gamma-1)\rho_0}{(\gamma+1)\rho_0-(\gamma-1)\rho_1},
\,\,\, \rho_1>\rho_0.
\end{equation}}
Given $\rho_0, p_0, \rho_1$, and $\gamma>1$, the other
variables $u_1$ and $p_1$ are determined by \eqref{2.10}. In
particular, the Mach number $M_1=\frac{u_1}{c_1}$ for state (1) is determined by
$M_1^2=\frac{2(\rho_1-\rho_0)^2}{\rho_0((\gamma+1)\rho_1-(\gamma-1)\rho_0)}$.

\begin{figure}[ht]
\centering
 \includegraphics[height=2.6in,width=2.2in]{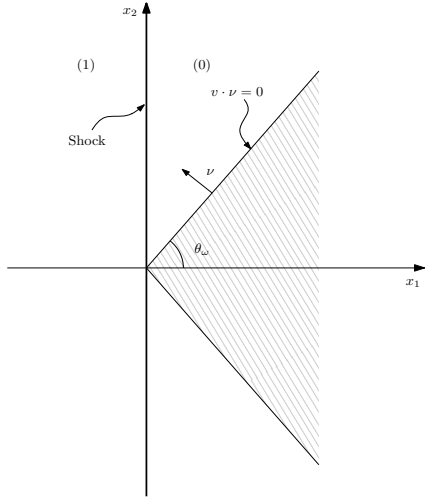}  
 \caption[]{Initial-boundary value problem}\protect\label{fig:IBV-1}
 \end{figure}

Since the initial-boundary value problem, {\it Problem 6.1},
is invariant under the self-similar scaling, we seek self-similar
solutions \eqref{self-similar}
governed by system \eqref{Euler1a}.
Since the problem is symmetric with respect to the axis $\eta=0$, it
suffices to consider the problem in the half-plane $\eta>0$ outside
the half-wedge:
$$
\Lambda:=\{\xi<0,\eta>0\}\cup\{\eta>\xi \tan\theta_w,\, \xi>0\}.
$$
Then {\it Problem 6.1} in the
$(t, {\bf x})$--coordinates can be formulated as the following
boundary value problem in the self-similar coordinates $(\xi,\eta)$:

\medskip
{\bf Problem 6.2 (Boundary value problem in the unbounded domain)}.
{\it Seek a solution to system \eqref{Euler1a} satisfying the
slip boundary condition on the wedge boundary:
$(U,V)\cdot\bnu =0$ on
$\partial\Lambda=\{\xi\le 0, \eta=0\}\cup\{\xi>0, \eta\ge
\xi\tan\theta_w\}$,
the asymptotic boundary condition as $\xi^2+\eta^2\to \infty$:
\begin{equation*}
{(\rho, U+\xi,V+\eta,p)\longrightarrow {\begin{cases}(\rho_0, 0,
0,p_0),\qquad
                         \xi>\xi_0, \eta>\xi\tan\theta_w,\\
              (\rho_1, u_1, 0, p_1),\qquad
                          \xi<\xi_0, \eta>0.
\end{cases}}
}
\end{equation*}}

\begin{figure}[h]
 \centering
 \includegraphics[height=1.9in,width=2.7in]{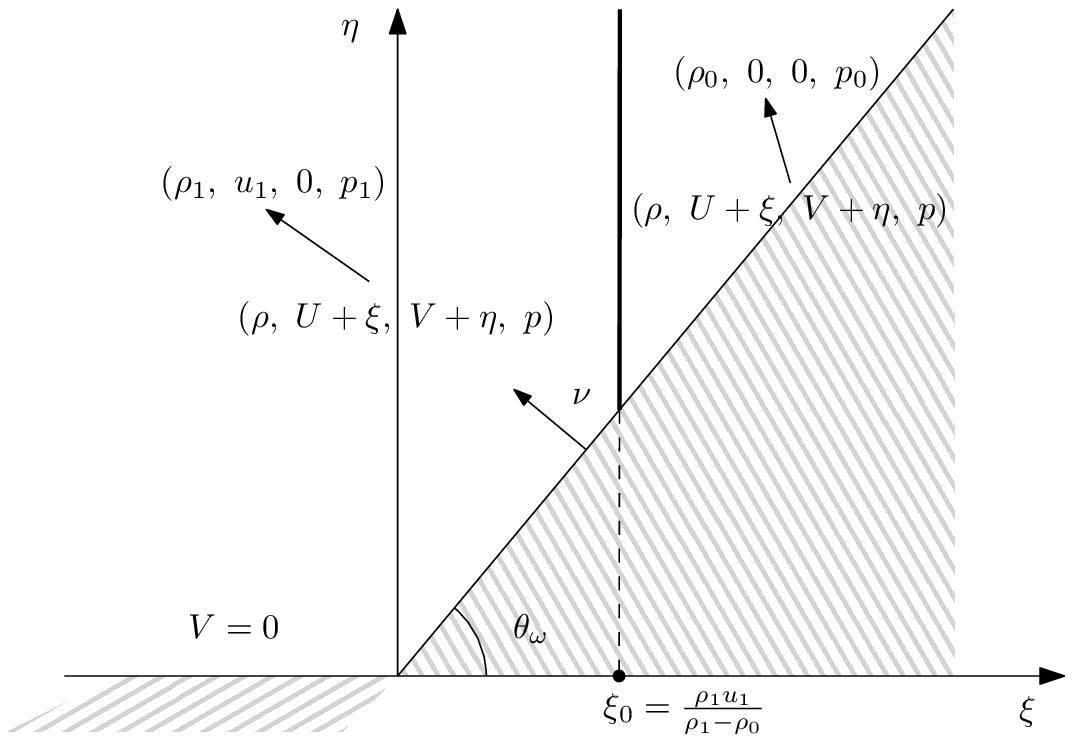}  
 \caption[]{Boundary value problem in the unbounded domain $\Lambda$}\label{fig:PM-1}
 \end{figure}
\begin{figure}[h]
 \centering
 \includegraphics[height=1.8in,width=2.2in]{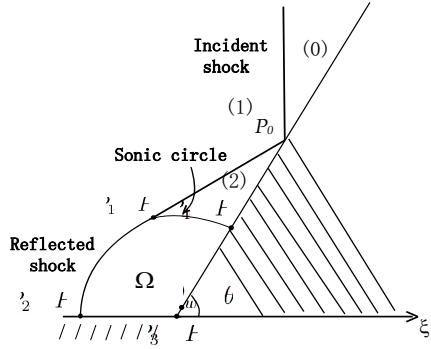}  
 \caption[]{Regular reflection-diffraction configuration}\label{fig:RC-1}
 \end{figure}

For our problem, since $\varphi_1$ does not satisfy the slip
boundary condition \eqref{boundary-condition-3}, the solution must
differ from $\varphi_1$ in $\{\xi<\xi_0\}\cap\Lambda$, thus a shock
diffraction-diffraction by the wedge vertex occurs.
The experimental, computational, and asymptotic analysis shows that various patterns of
reflected shock-fronts may occur, including the regular
and Mach reflections
(cf. \cite{BenDor,CFr,GM,Guderley,hunter1,HK,KB,KT,
Ligh,
Mach,Mora2,SA,TR,TH,TSK2,Timm,VanDyke,ZBHW}).
It is expected that the solutions of {\it Problem 6.2} contain all
possible patterns of shock reflection-diffraction configurations as
observed.
For a wedge angle $\theta_w\in (0, \frac{\pi}{2})$, different
reflection-diffraction patterns may occur. Various criteria and
conjectures have been proposed for the existence of configurations
for the patterns (cf. Ben-Dor \cite{BenDor}). One of the most important
conjectures made by von Neumann \cite{Neumann1,Neumann2} in 1943 is
the {\em detachment conjecture}, which states that the regular
reflection-diffraction configuration may exist globally whenever the
two-shock configuration (one is the incident shock-front and the other the
reflected shock-front) exists locally around the
point $P_0$ (see Fig. \ref{fig:RC-1}).
The following theorem was rigorously shown in Chang-Chen
\cite{CC} (also see Sheng-Yin \cite{ShengYin}, Bleakney-Taub
\cite{BT}, Neumann \cite{Neumann1,Neumann2}).

\smallskip
\begin{theorem}[Local theory]\label{local}
There exists $\theta_d=\theta_d(\rho_0,\rho_1,\gamma)\in (0, \frac{\pi}{2})$ such
that, when $\theta_w\in (\theta_d,\frac{\pi}{2})$, there are two states
$(2)$:$(\rho_2^a, U_2^a,V_2^a,p_2^a)$ and
$(\rho_2^b, U_2^b,V_2^b,p_2^b)$ such that
$|(U_2^a, V_2^a)|>|(U_2^b, V_2^b)|$ and
$|(U_2^b, V_2^b)| < c(\rho_2^b, S_2^b)$.
\end{theorem}

\smallskip
{\bf The von Neumann Detachment Conjecture} (\cite{Neumann1,Neumann2}):
{\em There exists a
global regular reflection-diffraction configuration whenever the
wedge angle $\theta_w$ is in $(\theta_d, \frac{\pi}{2})$}.

\smallskip
It is clear that the regular reflection-diffraction configuration is
not possible without a local two-shock configuration at the
reflection point on the wedge, so this is the weakest possible
criterion. In this case, the local theory indicates that there are
two possible choices for state (2). There had been a long debate to
determine which one is more physical for the local theory; see
Courant-Friedrichs \cite{CFr}, Ben-Dor \cite{BenDor}, and the references
cited therein.
Since the reflection-diffraction problem is not a local problem, we
take a different point of view that the selection of state (2)
should be determined by the global features of the problem, more
precisely, by the stability of the configuration with respect to the
wedge angle $\theta_w$, rather than the local features of the
problem.

\smallskip
{\bf Stability Criterion to Select the Correct
State (2) (Chen-Feldman \cite{CFe4})}:
{\it
Since the solution is unique when the wedge angle $\theta_w=\frac{\pi}{2}$,
it is required that the global regular reflection-diffraction
configuration be stable and converge to the unique normal
reflection solution when $\theta_w\to\frac{\pi}{2}$, provided that such a
global configuration can be constructed.}

\smallskip
We employ this stability criterion to conclude that our choice for
state (2) must be $(\rho_2^a, U_2^a, V_2^a, p_2^a)$. In general,
$(\rho_2^a, U_2^a, V_2^a, p_2^a)$ may be supersonic or subsonic. If
it is supersonic, the propagation speeds are finite and state (2) is
completely determined by the local information: state (1), state
(0), and the location of the point $P_0$. This is, any information
from the reflected region, especially the disturbance at the corner
$P_3$, cannot travel towards the reflection point $P_0$. However, if
it is subsonic, the information can reach $P_0$ and interact with
it, potentially altering the reflection-diffraction type. This
argument motivated the following  conjecture:

\smallskip
{\bf The von Neumann Sonic Conjecture \cite{Neumann1,Neumann2}}:
{\em There exists a
regular reflection-diffraction configuration when $\theta_w\in
(\theta_s, \frac{\pi}{2})$ for $\theta_s>\theta_d$ such that $|(U_2^a,
V_2^a)|>c_2^a$ at $P_0$.}

If state (2)
is sonic when $\theta_w=\theta_s$, then $|(U^a_2, V^a_2)|>c_2^a$ for
any $\theta_w\in (\theta_s, \frac{\pi}{2})$. This sonic conjecture is stronger
than the detachment one. In fact, the regime between the angles
$\theta_s$ and $\theta_d$ is very narrow and is only fraction of a
degree apart; see Sheng-Yin \cite{ShengYin}.

Following the argument in Chen-Feldman \cite{CFe5}, we have

\smallskip
\begin{theorem}[Chen-Feldman \cite{CFe5}]\label{p-dominate}
Let $(\rho,U,V,p)$ be a solution of {\it Problem 6.2} such that
$(\rho,U,V,p)$  is $C^{0,1}$ in the open region $P_0P_1P_2P_3$ and
the gradient of the tangential component of $(U,V)$ is continuous
across the sonic arc $\Sonic$.
Let $\Omega_1$ be the subregion of $\Omega$ formed by the fluid
trajectories past the sonic arc $\Gamma_{sonic}$. Then, in
$\Omega_1$, the potential flow equation
for self-similar solutions:
\begin{equation}\label{potential-1a}
\nabla\cdot \big(\rho(\nabla\varphi, \varphi)\nabla\varphi\big)
+2\rho(\nabla\varphi, \varphi)=0,
\end{equation}
with
$\rho(|\nabla\varphi|^2, \varphi)
=\big(\rho_0^{\gamma-1}-(\gamma-1)(\varphi+\frac{1}{2}|\nabla\varphi|^2)\big)^{\frac
1{\gamma-1}}$,
coincides with the full Euler equations
\eqref{Euler1a}, that is, equation \eqref{potential-1a}
is exact in the domain $\Omega_1$.
\end{theorem}

\smallskip
The regions such as $\Omega_1$ also exist in various Mach
reflection-diffraction configurations. Theorem \ref{p-dominate}
applies to such regions whenever the solution $(\rho,U,V,p)$ is
$C^{0,1}$ and the gradient of the tangential component of $(U,V)$ is
continuous. In fact, Theorem 6.2 indicates that, for the solutions
$\varphi$ of \eqref{potential-1a},
the $C^{1,1}$--regularity of $\varphi$ and the continuity of the tangential
component of the velocity field $(U,V)=\nabla\varphi$ are optimal
 across the sonic arc $\Gamma_{sonic}$.

Equation \eqref{potential-1a}
is a nonlinear
equation of mixed elliptic-hyperbolic type. It is elliptic if and
only if
$|\nabla\varphi| < c(|\nabla\varphi|^2,\varphi,\rho_0^{\gamma-1})$,
which is equivalent to
\begin{equation}
|\nabla \varphi| <q_*(\varphi, \rho_0, \gamma)
:=\sqrt{\frac{2}{\gamma+1}\big(\rho_0^{\gamma-1}-(\gamma-1)\varphi\big)}.
\label{1.1.8a}
\end{equation}

For the potential equation \eqref{potential-1a},
shock-fronts are discontinuities in the pseudo-velocity $\nabla\varphi$.
That is, if $D^+$ and $D^-:=D\setminus\overline{D^+}$ are two
nonempty open subsets of $D\subset\R^2$, and $\mathcal{S}:=\partial D^+\cap D$
is a $C^1$--curve where $D\varphi$ has a jump, then $\varphi\in
W^{1,1}_{loc}(D)\cap C^1(D^\pm\cup S)\cap C^2(D^\pm)$ is a global
weak solution of (\ref{potential-1a})
in $D$ if
and only if $\varphi$ is in $W^{1,\infty}_{loc}(D)$ and satisfies
equation \eqref{potential-1a} in $D^\pm$ and the Rankine-Hugoniot
condition on $\mathcal{S}$:
\begin{equation}\label{FBConditionSelfSim-0}
\left[\rho(|\nabla\varphi|^2,\varphi)\nabla\varphi\cdot\bnu\right]_{\mathcal S}=0.
\end{equation}

Then the plane incident shock-front solution in the $(t, {\bf
x})$--coordinates with states $(\rho, \nabla_{\bf x}\Phi)=(\rho_0,
0,0)$ and $(\rho_1, u_1,0)$ corresponds to a continuous weak
solution $\varphi$ of (\ref{potential-1a}) in the self-similar
coordinates $(\xi,\eta)$ with the following form:
\begin{eqnarray}
&&\varphi_0(\xi,\eta)=-\frac{1}{2}(\xi^2+\eta^2) \qquad
 \hbox{for } \,\, \xi>\xi_0,
 \label{flatOrthSelfSimShock1} \\
&&\varphi_1(\xi,\eta)=-\frac{1}{2}(\xi^2+\eta^2)+ u_1(\xi-\xi_0)
\qquad
 \hbox{for } \,\, \xi<\xi_0,
 \label{flatOrthSelfSimShock2}
\end{eqnarray}
respectively, where
$u_1
=\sqrt{\frac{2(\rho_1-\rho_0)(\rho_1^{\gamma-1}-\rho_0^{\gamma-1})}
 {(\gamma-1)(\rho_1+\rho_0)}}>0$ and
$\xi_0
=\frac{\rho_1u_1}{\rho_1-\rho_0}>0
$
are the velocity of state (1) and the location of the incident
shock-front, uniquely determined by $(\rho_0,\rho_1,\gamma)$ through
(\ref{FBConditionSelfSim-0}). Then
 $P_0=(\xi_0, \xi_0 \tan\theta_w)$ in Fig. 2, and
{\it Problem 6.2} in the context of the potential flow equation can
be formulated as

\medskip
{\bf Problem 6.3}\label{BVP} {(Boundary value problem)} (cf. Fig.
{\rm 2}). {\it Seek a solution $\varphi$ of
\eqref{potential-1a}
in the self-similar domain
$\Lambda$ with the boundary condition on
$\partial\Lambda$:
\begin{equation}\label{boundary-condition-3}
\nabla\varphi\cdot\bnu|_{\partial\Lambda}=0,
\end{equation}
and the asymptotic boundary condition at infinity:
\begin{equation}\label{boundary-condition-2}
\varphi\to\bar{\varphi}:=
\begin{cases} \varphi_0 \qquad\mbox{for}\,\,\,
                         \xi>\xi_0, \eta>\xi \tan\theta_w,\\
              \varphi_1 \qquad \mbox{for}\,\,\,
                          \xi<\xi_0, \eta>0,
\end{cases}
\qquad \mbox{when $\xi^2+\eta^2\to \infty$},
\end{equation}
where {\rm (\ref{boundary-condition-2})} holds in the sense that
$
\displaystyle
\lim_{R\to\infty}\|\varphi-\overline{\varphi}\|_{C(\Lambda\setminus
B_R(0))}=0.
$
}

\smallskip
 In Chen-Feldman
\cite{CFe4}, we first followed the von
Neumann criterion and the stability criterion introduced above
to establish a local existence theory of regular shock reflection
near the reflection point $P_0$ in the level of potential flow,
when the wedge angle is large and close to $\frac{\pi}{2}$.  In this case,
the vertical line is the incident shock-front
that hits
the wedge at $P_0=(\xi_0, \xi_0 \tan\theta_w)$, and state
(0) and state (1) ahead of and behind $S$ are given by $\varphi_0$
and $\varphi_1$ defined in \eqref{flatOrthSelfSimShock1} and
\eqref{flatOrthSelfSimShock2}, respectively. The solutions $\varphi$
and $\varphi_1$ differ only in the domain $P_0\PtUpL \PtLwL \PtLwR$
because of shock diffraction by the wedge, where the curve
$P_0\PtUpL \PtLwL$ is the reflected shock-front with the straight segment
$P_0\PtUpL$. State (2) behind $P_0\PtUpL$ can be computed explicitly
with the form:
\begin{equation}\label{state2}
\varphi_2(\xi,\eta)=-\frac{1}{2}(\xi^2+\eta^2)+u_2(\xi-\xi_0)+
(\eta-\xi_0\tan\theta_w)u_2\tan\theta_w,
\end{equation}
which satisfies $ \nabla\varphi\cdot \bnu=0$ on $\partial\Lambda\cap
\{\xi>0\}$;  the constant velocity $u_2$ and the angle $\theta_s$
between $P_0\PtUpL$ and the $\xi$--axis are determined by
$(\theta_w,\rho_0,\rho_1,\gamma)$ from the two algebraic equations
expressing (\ref{FBConditionSelfSim-0}) and continuous  matching of
states (1) and (2) across $P_0\PtUpL$
as in Theorem 6.1.
Moreover, $\varphi_2$ is the unique solution in the domain
$P_0\PtUpL \PtUpR$, as argued in \cite{CC,Serreref}. Hence
$
P_1P_4:=\Sonic=\partial\Omega\cap\partial B_{c_2}(u_2,
u_2\tan\theta_w)
$
is the sonic arc of state $(2)$ with center $(u_2,
u_2\tan\theta_w)$ and radius $c_2$.

It should be noted that, in order that the solution $\varphi$ in the
domain $\Omega$ is a part of the global solution to {\it Problem
6.3}, that is, $\varphi$ satisfies the equation in the sense of
distributions in $\Lambda$, especially across the sonic arc
$P_1P_4$, it is required that
$
\nabla(\varphi-\varphi_2)\cdot \bnu|_{P_1P_4}=0.
$
That is, we have to match our solution with state (2), which is the
necessary condition for our solution in the domain $\Omega$ to be a
part of the global solution. To achieve this, we have to show that
our solution is at least $C^1$ with $\nabla (\varphi-\varphi_2)=0$
across $P_1P_4$.
Then {\it Problem 6.3} can be reformulated as the following free boundary
problem:

\medskip
{\bf Problem 6.4 (Free boundary problem)}. {\it Seek a solution
$\varphi$ and a free boundary $P_1P_2=\{\xi=f(\eta)\}$ such that
\begin{itemize}
\item[\rm (i)] $f\in C^{1,\alpha}$ and
\begin{equation}\label{7.11}
\Omega:=\{\xi>f(\eta)\}\cap \Lambda{=\{\varphi<\varphi_1\}\cap \Lambda};
\end{equation}

\item[\rm (ii)] $\varphi$ satisfies the free boundary condition \eqref{FBConditionSelfSim-0}
along $P_1P_2$;

\item[\rm (iii)] $\varphi\in C^{1,\alpha}(\overline{\Omega})\cap C^2(\Omega)$
     solves \eqref{potential-1a} in $\Omega$, is subsonic in $\Omega$,
     and satisfies
\begin{eqnarray}
&&(\varphi-\varphi_2,
\nabla(\varphi-\varphi_2)\cdot\bnu)|_{P_1P_2}=0,
  \label{7.12}\\
&&\nabla \varphi\cdot\bnu|_{P_3P_4\cup\Gamma_{symm}}=0. \label{7.13}
\end{eqnarray}
\end{itemize}
}

The boundary condition on $\Gamma_{symm}$ implies that $f'(0)=0$ and
thus ensures the orthogonality of the free boundary with the
$\xi$-axis. Formulation \eqref{7.11} implies that the free boundary
is determined by the level set $\varphi=\varphi_1$.
The
free boundary condition \eqref{FBConditionSelfSim-0} along $P_1P_2$ is the conormal
boundary condition on $P_1P_2$. Condition \eqref{7.12} ensures that
the solution of the free boundary problem in $\Omega$ is a part of
the global solution. Condition \eqref{7.13}
is the slip boundary condition.
{\it Problem 6.4} involves two types of transonic flow: one is a
continuous transition through the sonic arc $P_1P_4$ as a fixed
boundary from the supersonic region (2) to the
subsonic region $\Omega$; the other is a jump transition
through the transonic shock-front as a free boundary from the supersonic
region (1) to the subsonic region $\Omega$.

In Chen-Feldman \cite{CFe4,CFe5,CFe6},
we have
developed a rigorous mathematical approach to solve the von Neumann
sonic conjecture, {\it Problem
6.4}, and established a global theory for solutions of regular
reflection-diffraction up to the sonic angle,
which converge to the unique solution of the
normal shock reflection when $\theta_w$ tends to ${\pi}/{2}$.
For more details, see Chen-Feldman \cite{CFe4,CFe5,CFe6}.

The mathematical existence of Mach reflection-diffraction
configurations is still open; see \cite{BenDor,ChS6,CFr}.
Some progress has been made in the recent years in the study
of the 2-D Riemann problem for hyperbolic conservation laws;
see \cite{CCY2,CH,Chen1,GJLZ,GK,GM,GSZ,KTa,LL,LZY,LZheng,SCG,Ser3,Serre-09,SZ,TZ,Zhe,Zhe1}
and the references cited therein.

Some recent developments in the study of the local nonlinear
stability of multidimensional compressible vortex sheets can be found in
Coulombel-Secchi \cite{CS}, Bolkhin-Trakhinin \cite{BT},
Trakhinin \cite{Tra},
Chen-Wang \cite{CWangYG}, and the references cited therein.
Also see Artola-Majda \cite{AM}.
For the construction of the non-self-similar global solutions
for some multidimensional systems, see Chen-Wang-Yang \cite{ChenWangYang} and
the references cited therein.

\section{Divergence-Measure Vector Fields and Multidimensional Conservation Laws}
Naturally, we want to approach the questions of existence, stability,
uniqueness, and long-time behavior of entropy solutions for
M-D hyperbolic conservation laws
{with neither specific reference to any particular method for
constructing the solutions nor additional regularity assumptions}.
Some recent efforts have been in developing
a theory of divergence-measure fields to construct a global framework
for the analysis of solutions of M-D
hyperbolic systems of conservation
laws.

Consider  system \eqref{cons}
in $\R^d$.
As discussed in \S 2.5,
the $BV$ bound generically fails for the M-D case.
In general, for M-D conservation laws, especially
the Euler equations,
solutions of \eqref{cons} are expected to be in the following
class of entropy solutions:

\medskip
\begin{enumerate}
\item[(i)]  $\uu(t,\xx)\in \M(\R_+^{d+1})$,
   or $L^p(\R_+^{d+1}),\,  1\le p\le \infty$;

\smallskip
\item[(ii)] $\uu(t,\xx)$ satisfies the Lax entropy inequality:
\begin{equation}\label{entropy-ineq}
\mu_{\eta}:=\partial_t\eta(\uu(t,\xx))
+\nabla_\xx\cdot\qq(\u(t,\xx)) \le 0
\end{equation}
in the sense of distributions for any convex entropy pair
$(\eta,\qq): \R^n\to \R\times \R^d$
so that  $\eta(\uu(t,\xx))$ and $\qq(\uu(t,\xx))$ are distributions.
\end{enumerate}

\noindent
The Schwartz Lemma infers from (\ref{entropy-ineq})
that the distribution $\mu_\eta$ is in fact a Radon measure:
$
\hbox{div}_{(t,\xx)}(\eta(\uu(t,\xx)),\qq(\uu(t,\xx)))
\in \M(\R_+^{d+1}).
$
Furthermore, when $\uu\in L^\infty$, this is also true for any $C^2$
entropy-entropy flux pair $(\eta,\qq)$ ($\eta$ not necessarily
convex) if (\ref{cons}) has a strictly convex entropy,
which was first observed in Chen \cite{Ch3}.
More generally, we have

\medskip
\noindent
{\bf Definition.} Let $\D\subset{\R}^N$ be open.
For $1\le p\le\infty$, $\FF$ is called a $\DM^p(\D)$--field if
$\FF\in L^p(\D;{\R}^N)$ and
\begin{equation}\label{norm1}
\|\FF\|_{\DM^p(\D)}
:=\|\FF\|_{L^p(\D;\R^N)}+\|\hbox{div}\FF\|_{\M(\D)}<\infty;
\end{equation}
and the field $\FF$ is called a ${\DM}^{ext}(\D)$--field
if $\FF\in\M(\D;\R^N)$
and
\begin{equation}\label{norm2}
\|\FF\|_{\DM^{ext}(\D)}
:=\|(\FF, \hbox{div} \FF)\|_{\M(\D)}<\infty.
\end{equation}
Furthermore, for any bounded open set $\D\subset{\R}^N$,
$\FF$ is called a ${\DM}^p_{loc}({\R}^N)$--field
if $\FF\in {\DM}^p(\D)$;
and $\FF$ is called a ${\DM}^{ext}_{loc}({\R}^N)$--field
if $\FF\in {\DM}^{ext}(\D)$.
A field $\FF$ is simply called a $\DM$--field in $\D$ if
$\FF\in \DM^p(\D), 1\le p\le \infty$, or $\FF\in\DM^{ext}(\D)$.

\smallskip
It is easy to check that these spaces, under the respective norms
$\|\FF\|_{\DM^p(\D)}$ and $\|\FF\|_{\DM^{ext}(\D)}$, are Banach spaces.
These spaces are larger than the space of $BV$-fields.
The establishment of the Gauss-Green theorem, traces, and other
properties of $BV$ functions in the 1950s
(cf. Federer \cite{Fe}; also Ambrosio-Fusco-Pallara \cite{A1},
Giusti \cite{G1}, and Volpert \cite{Vo}) has significantly advanced
our understanding of solutions of nonlinear PDEs and related problems
in the calculus of variations,
differential geometry, and other areas,
especially for the 1-D theory of hyperbolic conservation laws.
A natural question is whether the $\DM$--fields have similar
properties, especially the normal traces and the Gauss-Green formula
to deal with entropy solutions for M-D conservation laws.

On the other hand,
motivated by various nonlinear problems from conservation laws,
as well as for rigorous derivation of systems of balance laws with
measure source terms from the physical principle of balance law
and the recovery of Cauchy entropy flux through
the Lax entropy inequality for entropy solutions of hyperbolic
conservation laws by capturing the entropy dissipation,
a suitable notion of normal traces and corresponding Gauss-Green formula
for divergence-measure fields are required.

Some earlier efforts were made on generalizing the Gauss-Green
theorem for some special situations, and relevant results can be found in
Anzellotti \cite{An} for an abstract formulation for $\FF\in L^\infty$,
Rodrigues \cite{Ro} for $\FF\in  L^2$, and
Ziemer \cite{Zi1} for a related problem for $\hbox{div}\FF\in L^1$;
also see Baiocchi-Capelo \cite{BC} and Brezzi-Fortin \cite{BF}.
In Chen-Frid \cite{CF6}, an explicit way to calculate the suitable
normal traces was first observed for $\FF\in \DM^\infty$,
under which a generalized Gauss-Green theorem was shown to hold,
which has motivated the development of
a theory of divergence-measure fields
in \cite{CF6,ChenTorres,ChenTorresZiemer1}.

Some entropy methods based on the theory of divergence-measure
fields presented above have been developed and applied for solving
various nonlinear problems
for hyperbolic conservation laws and related nonlinear PDEs.
These problems especially include
(i) Stability of Riemann solutions, which may
contain rarefaction waves, contact discontinuities, and/or vacuum states,
in the class of entropy solutions of the Euler equations
in \cite{Chen-Chen,CF6,CFL};
(ii) Decay of periodic entropy solutions
in \cite{CF2};
(iii) Initial and boundary layer problems
  in \cite{CF6,CR,ChenTorres,Vas};
(iv) Rigorous derivation of systems of balance laws
   from the physical principle of balance law
  and the recovery of Cauchy entropy flux through
  the Lax entropy inequality for entropy solutions of hyperbolic
  conservation laws by capturing the entropy dissipation
    in \cite{ChenTorresZiemer1};

\smallskip
It would be interesting to develop further the theory of divergence-measure fields
and more efficient
entropy methods for solving more various problems in PDEs and related areas
whose solutions are only measures or $L^p$ functions.
For more details, see \cite{CF6,ChenTorres,ChenTorresZiemer1}.

\smallskip
We also refer the reader to the related papers on other aspects of multidimensional
conservation laws in this volume.

\bigskip
\bigskip
{\bf Acknowledgements.}  The work of Gui-Qiang G. Chen was supported in part by NSF grants
DMS-0935967, DMS-0807551, the Royal Society--Wolfson Research
Merit Award (UK), and the EPSRC Science and Innovation
award to the Oxford Centre for Nonlinear PDE (EP/E035027/1).

\end{document}